\documentclass[a4paper,10pt]{article}

\usepackage{amsmath,amssymb,amsthm,xypic,epsfig}

\def\g{\mathfrak{g}}
\def\h{\mathfrak{h}}
\def\n{\mathfrak{n}}
\def\b{\mathfrak{b}}
\def\C{\mathbb{C}}

\def\Z{\mathbb{Z}}
\def\qed{$\hfill \blacksquare$}
\def\l{\mathfrak{l}}

\newtheorem{theo}{Theorem}[section]
\newtheorem{prop}{Proposition}[section]
\newtheorem{lem}{Lemma}[section]
\newtheorem{cor}{Corollary}[section]
\numberwithin{equation}{section}

\begin{document}
 \title{Twisted Traces of Quantum Intertwiners and
Quantum Dynamical R-Matrices Corresponding to Generalized
Belavin-Drinfeld Triples}
 \author{P. Etingof and O. Schiffmann}
\date{}
 \maketitle
\section{Introduction}
 This paper is a continuation of \cite{ES1} and \cite{EV4}. 

  In \cite{EV4}, A.Varchenko and the first author considered 
weighted traces of products of intertwining operators for quantum groups 
$U_q(\g)$, where $\g$ is a simple Lie algebra. 
They showed that the generating function 
$F_{V_1,\ldots V_N}(\lambda,\mu)$  of such traces (where $\lambda,\mu$ are
complex weights for $\g$) 
satisfies four commuting systems of difference equations -- 
the Macdonald-Ruijsenaars (MR) system, the quantum 
Knizhnik-Zamolodchikov-Bernard (qKZB) system, the dual MR
system, and the dual qKZB system. The first two systems 
are systems of difference equations with respect to $\lambda$,
which involve Felder's trigonometric dynamical R-matrix 
depending of $\lambda$.
The second two systems are systems of difference equations respect to $\mu$, 
which are obtained from the first two by 
the transformation $\lambda\to \mu,V_i\to V_{N-i+1}^*$. 
Such a symmetry is explained by the fact 
that the function $F_{V_1,\ldots V_N}(\lambda,\mu)$ is invariant 
under this transformation. 

  If the quantum group $U_q(\g)$ is replaced with the Lie algebra $\g$, 
these results are replaced with their classical analogs (\cite{EV4}). 
Namely, the MR and qKZB equations are replaced by the 
classical MR and KZB equations, which are differential equations 
involving Felder's classical trigonometric 
dynamical r-matrix. The dual MR and KZB 
equations retain roughly the same form, but involve the 
rational quantum dynamical 
R-matrix rather than the trigonometric one. Thus, the symmetry between 
$\lambda$ and $\mu$ is destroyed. 

  In \cite{ES1}, we generalized the classical MR and KZB equations 
to the case when the trace is twisted using a
"generalized Belavin-Drinfeld triple", i.e. 
a pair of subdiagrams $\Gamma_1,\Gamma_2$ of the Dynkin diagram 
of $\g$ together with an isomorphism $T:\Gamma_1\to \Gamma_2$ between them.
It turned out that such twisted traces also satisfy 
differential equations which involve a dynamical r-matrix, namely
the one attached to the triple $(\Gamma_1,\Gamma_2,T)$ by the 
second author in \cite{S3}. 

  After \cite{ES1} was finished, we wanted to generalize its results 
to the quantum case. It was clear to us that to express the result
we would need an explicit quantization of classical dynamical r-matrices   
from \cite{S3}. Therefore, we hoped that attempts to quantize 
the results of \cite{ES1} using the approach of \cite{EV4} 
could help us obtain such a quantization (which was unknown even for the
usual Belavin-Drinfeld classical 
r-matrices). This program did, in fact, succeed,
and the quantization of dynamical r-matrices from \cite{S3}
was recently obtained in \cite{ESS}. 

  In this paper, using the results of \cite{ESS} and methods of \cite{EV4}, we generalize 
the difference equations from \cite{EV4} to the twisted case; this 
also provides a quantum generalization 
of \cite{ES1}. Namely, we deduce difference equations with respect to the 
weight 
$\lambda$    
for the generating function of the twisted traces for $U_q(\g)$, -
the twisted MR and qKZB equations. Not surprisingly, they involve 
the dynamical R-matrix constructed in \cite{ESS}. In the case when $T$ 
is an automorphism of the whole Dynkin diagram of $\g$, 
we also deduce the twisted dual MR and qKZB equations, i.e. the
difference equations with respect to the other weight $\mu$. 
These equations involve the usual (Felder's) dynamical R-matrix, 
but differ from the equations of \cite{EV4} by explicit occurrence of $T$. 
Thus, we see that for $T\ne 1$, there is no symmetry between $\lambda$ and 
$\mu$. 

  If $T$ is not an automorphism, we do not expect the existence of the dual 
equations. This is explained at the end of Section 2.
 
  Replacing $U_q(\g)$ with $\g$, we obtain the classical limit
of these results. The twisted MR and qKZB equations become their classical 
analogs from \cite{ES1}. The dual equations retain their form, but 
the trigonometric R-matrices are replaced by their rational limits. 

  Finally, we adapt the construction of the quantum dynamical
R-matrices from \cite{ESS} to the case when $\g$ is an arbitrary symmetrizable
Kac-Moody algebra. This yields quantizations of the classical dynamical 
r-matrices from \cite{ES1} in the case of Kac-Moody algebras. Again, 
the generating functions for twisted traces of intertwiners for $U_q(\g)$
satisfy a set of difference equations involving these quantum dynamical
R-matrices, and a set of dual difference equations if in addition $T$ is
an automorphism of the Dynkin diagram.

  In the next paper, we plan to generalize these results 
to the case of affine algebras, when traces take values in finite-dimensional 
representations. This involves dynamical R-matrices with \textit{spectral 
parameters}. In particular, we plan to obtain 
a trace representation of solutions of the elliptic qKZ equation
(with Belavin's elliptic R-matrix), and compute its monodromy. 

\paragraph{Remark.} The elliptic qKZ
 equation is important in statistical mechanics (see \cite{JM}).
For its classical version, the trace representation of solutions 
and monodromy are obtained in \cite{E1}. The problem of quantizing the results
of \cite{E1} was suggested to the first author by his advisor I. Frenkel as a topic for his PhD thesis in 1992. After this the first author
tried to quantize the results of \cite{E1} (see \cite{E2}) but 
obtained only partial results.
 \subsection{Notations}
 \paragraph{}Let $\g$ be a simple complex Lie algebra with a fixed polarization
 $\g=\n_- \oplus \h
\oplus \n_+$. Let $\Gamma$
(resp. $\Delta$) be the Dynkin diagram (resp. the root system) of
$\g$. Denote by $(a_{ij})$ the Cartan matrix of $\g$ and let $d_i$ be
relatively prime positive integers such that $(d_ia_{ij})$ is a
symmetric matrix. Let $(\,,\,)$ be the nondegenerate invariant symmetric form
for which  $(\alpha,\alpha)=2$ if $\alpha$ is a \textit{short} root.
 Let $\{e_\alpha,f_\alpha\}_{\alpha \in \Delta}$ be a
Chevalley basis of $\n_- \oplus \n_+$, normalized in such a way that
$(e_\alpha,f_\alpha)=1$ for all $\alpha$ and set
$h_\alpha=[e_\alpha,f_\alpha]$. We also let $\Omega \in \g \otimes \g$
and $\Omega_\h \in \h \otimes \h$ be the inverse elements of the restriction of $(\,\,)$ to $\g$ and
$\h$ respectively.\\ 
\paragraph{}Let $q=e^{\frac{1}{2}\hbar}$ where $\hbar$ is a formal
variable. For any operator $A$ we set
$q^A=e^{\hbar\frac{A}{2}}$. Let $U_{q}(\g)$ be the Drinfeld-Jimbo
quantized enveloping algebra of $\g$. It is a $\C[[\hbar]]$-Hopf algebra
with generators $E_\alpha,F_\alpha$, $\alpha \in \Gamma$ and $q^h$, $h
\in \h$ subject to the following set of relations :
 $$q^{x+y}=q^xq^y,\
x,y\in \h\quad q^hE_{\alpha_j}q^{-h}=q^{\alpha_j(h)}E_{\alpha_j},\quad
q^hF_{\alpha_j}q^{-h}=q^{-\alpha_j(h)}F_{\alpha_j}$$
$$E_{\alpha_i}F_{\alpha_j}-F_{\alpha_j}E_{\alpha_i}=\delta_{ij}\frac{q^{d_ih_{\alpha_i}}-q^{-d_ih_{\alpha_i}}}{q^{d_i}-q^{-d_i}},$$
$$\sum_{k=0}^{1-a_{ij}} (-1)^k \bmatrix 1-a_{ij} \\ k
\endbmatrix_{q^{d_i}} E_{\alpha_i}^{1-a_{ij}-k}E_{\alpha_j}E_{\alpha_i}^k=0,\quad i\ne j,$$
$$\sum_{k=0}^{1-a_{ij}} (-1)^k \bmatrix 1-a_{ij} \\ k
\endbmatrix_{q^{d_i}} F_{\alpha_i}^{1-a_{ij}-k}F_{\alpha_j}F_{\alpha_i}^k=0,\quad i\ne j.$$
 where as usual
 $$ \bmatrix n \\ k \endbmatrix
_q=\frac{[n]_q!}{[k]_q! [n-k]_q!}, \quad [n]_q! = [1]_q \cdot [2]_q
\cdot \dots \cdot [n]_q, \quad [n]_q= \frac {q^n - q^{-n}}{q-q^{-1}}$$
\paragraph{}Comultiplication $\Delta,$ antipode $S$ and counit
$\epsilon$ in $U_{q}(\g)$ are given by 
$$\Delta(E_{\alpha_i}) = E_{\alpha_i}\otimes
q^{d_ih_{\alpha_i}} + 1\otimes E_{\alpha_i}, \quad
  \Delta(F_{\alpha_i}) = F_{\alpha_i}\otimes 1 + q^{-d_ih_{\alpha_i}}\otimes F_{\alpha_i}, \quad
  \Delta(q^h) = q^h \otimes q^h$$
 $$S(E_{\alpha_i})=-E_{\alpha_i}q^{-d_ih_{\alpha_i}},\quad
S(F_{\alpha_i})=-q^{d_ih_{\alpha_i}}F_i,\quad S(q^h)=q^{-h}$$
 $$\epsilon(E_{\alpha_i}) =
\epsilon(F_{\alpha_i}) = 0,\quad \epsilon(q^h) = 1.$$ 
\paragraph{}Let
$U_{q}(\n_{\pm})$ be the subalgebra generated by $(E_\alpha)_{\alpha\in \Gamma}$ and
$(F_\alpha)_{\alpha\in \Gamma}$ respectively. It is known that $U_{q}(\g)$ is
quasitriangular, with R-matrix $\mathcal{R} \in q^{\Omega_\h} U_{q}(\n_+)
\hat{\otimes} U_{q}(\n_-)$. Here $\hat{\otimes}$ denotes the completion
with respect to the principal grading of $U_{q}(\n_{\pm})$.

\subsection{Generalized Belavin-Drinfeld triples and classical dynamical 
r-matrices}
\paragraph{}Let $\l \subset \h$ be a subalgebra on which $(\,,\,)$ is 
nondegenerate. Let $(x_i)_{i \in I}$ be an orthonormal basis of $\l$ and let
 $(x^i)_{i\in I}$ be the dual basis of $\l^*$. The classical dynamical Yang-Baxter 
equation with respect to $\l$ is the following equation :
\begin{equation}\label{E:CDYBE}
\begin{split}
\sum_i &\left(x_i^{(1)} \frac{\partial r^{23}(\lambda)}{\partial x^i}-
x_i^{(2)} \frac{\partial r^{13}(\lambda)}{\partial x^i} + x_i^{(3)} 
\frac{\partial r^{12}(\lambda)}{\partial x^i}\right)\\
& +
[r^{12}(\lambda),r^{13}(\lambda)]+[r^{12}(\lambda),r^{23}(\lambda)]+
[r^{13}(\lambda),r^{23}(\lambda)]=0
\end{split}
\end{equation}
where $r(\lambda): \l^* \to (\g \otimes \g)^\l$ is a meromorphic function.
 Solutions of (\ref{E:CDYBE}) relevant to the theory of Poisson-Lie groupoids
 (see \cite{EV1}, \cite{ES6}, \cite{Xu}) are those satisfying the generalized unitarity 
condition, i.e $r(\lambda) + r^{21}(\lambda)=\Xi$ is constant and $\Xi$ belongs to $(S^2\g)^\g$. In \cite{S3} the second author classified all
 such solutions $r(\lambda)$ which are non skewsymmetric (that is, $\Xi \neq 0$).
 Up to isomorphism and 
gauge 
transformations, they are labeled by the following combinatorial data called generalized Belavin-Drinfeld triples.
\paragraph{Definition.} A generalized Belavin-Drinfeld triple is a triple 
$(\Gamma_1,\Gamma_2,T)$ where $\Gamma_1,\Gamma_2 \subset \Gamma$ and $T:
 \Gamma_1 \stackrel{\sim}{\to} \Gamma_2$ is an orthogonal isomorphism.
\paragraph{}Let $(\Gamma_1,\Gamma_2,T)$ be a generalized 
Belavin-Drinfeld triple. Set
$$\l=\big(\sum_\alpha \C(\alpha -T(\alpha))\big)^\perp \subset \h.$$
Note that $\l$ is spanned by real elements so that the restriction of $(\,,\,)$
 to $\l$ is nondegenerate. Let $\h_0\subset \h$ be the orthogonal complement
 to $\l$ in $\h$ and let $\Omega_{\h_0} \in \h_0 \otimes \h_0$ be the
element inverse to the form $(\,,\,)$. The following lemmas are proved in \cite{ES1}.

\begin{lem}\label{L:00001} There exists a unique Lie algebra homomorphism $B:\b_-
\to\b_-$ (resp.
$B^{-1}: \b_+ \to \b_+$) such that
$B(f_\alpha)=f_{T(\alpha)}$, $B(h_\alpha)=h_{T(\alpha)}$ if $\alpha \in 
\Gamma_1$,
$B(f_\alpha)=0$ if
$\alpha
\not\in \Gamma_1$ (resp.
$B^{-1}(e_\alpha)=e_{T^{-1}(\alpha)}$, $B^{-1}(h_\alpha)=h_{T^{-1}(\alpha)}$
 if
$\alpha \in \Gamma_2$, $B^{-1}(e_\alpha)=0$ if
$\alpha
\not\in \Gamma_2$), and $B^{\pm 1}(h)=h$ if $h \in \l$. Moreover the restriction of $B$ to $\h$ is an orthogonal operator.\end{lem}
\paragraph{Remark.} We use the symbol $B^{-1}$ for notational convenience 
only. The operators $B$ and $B^{-1}$ are only inverse to each other when 
restricted to $\h$.

\begin{lem}[Cayley transform] For any $x\in \h_0$, there exists a unique 
element $C_T(x)\in \h_0$ such that 
for all $\alpha\in \Gamma_1$ one has $(\alpha-T\alpha,C_T(x))=(\alpha+
T\alpha,x)$. The linear operator $C_T: \h_0 \to \h_0$ is skew-symmetric.
\end{lem}
The classical dynamical r-matrix associated to $(\Gamma_1,\Gamma_2,T)$ is
\begin{equation}
r_T(\lambda)=-r_0^{21}+\sum_{\alpha,l\geq 1} e^{-l(\alpha,\lambda)}e_\alpha \wedge B^l f_\alpha +\frac{1}{2}(C_T \otimes 1)\Omega_{\h_0}
\end{equation}
where $r_0=\frac{1}{2}\Omega_\h + \sum_\alpha e_\alpha \otimes f_\alpha$ is the standard classical r-matrix.
\subsection{Quantum dynamical R-matrices}
\paragraph{}In our joint work with Travis Schedler \cite{ESS} we obtain an explicit quantization of
 the r-matrices $r_T(\lambda)$. Namely, we construct a trigonometric
function
 $\widetilde{R}_T(\lambda): \l^* \to U_{q}(\g) {\otimes} U_{q}(\g)$ (tensor product in the category of topologically free $\C[[\hbar]]$-modules) such that
 $\widetilde{R}_T(\lambda)\equiv 1-\hbar r_T(\lambda)\;\mathrm{mod}\;(\hbar^2)$
 which satisfies the quantum dynamical Yang-Baxter equation
\begin{equation}\label{E:QDYBE}
\begin{split}
\widetilde{R}_T^{12}(\lambda-\frac{1}{2}\hbar h^{(3)})\widetilde{R}_T^{13}
(\lambda&+\frac{1}{2}\hbar h^{(2)})\widetilde{R}_T^{23}(\lambda-\frac{1}{2}
\hbar h^{(1)})\\
&=\widetilde{R}_T^{23}(\lambda+\frac{1}{2}\hbar h^{(1)})\widetilde{R}_T^{13}
(\lambda-\frac{1}{2}\hbar h^{(2)})\widetilde{R}_T^{12}(\lambda+\frac{1}{2}\hbar
 h^{(3)}).
\end{split}
\end{equation}
In the above equation we used the usual notation for shifts in the dynamical variable: for 
instance, if $S(\lambda)$ is any meromorphic function $\l^* \to 
U_{q}(\g)^{\otimes 2}$ we set $S(\lambda-\frac{1}{2}\hbar h^{(3)})=S(\lambda)-
\frac{1}{2}\hbar \sum_i \frac{\partial S}{\partial y^i} y_i^{(3)}+\ldots$ (the
 Taylor expansion), where $y_1,\ldots,y_r$ is a basis of $\l$ and $y^1,\ldots
 y^r$ is the dual basis of $\l^*$.
\paragraph{}The construction is based on the following result. Let $I_\pm 
\subset U_{q}(\b_\pm)$ be the kernels of the projections $U_{q}(\b_\pm) \to
 U_{q}(\h)$. Also set $Z=\frac{1}{2}((1-C_T) \otimes 1) \Omega_{\h_0}$. The maps $B: U_{q}(\b_-)\to U_{q}(\b_-)$ and $B^{-1}: U_{q}(\b_+)\to U_{q}(\b_+)$ are defined in the same fashion as in the classical case (see 
Lemma~\ref{L:00001}). 
\paragraph{}To simplify notations we will write $q_i^{A}$ for $(q^A)_i$ for any operator $A$ (the operator $q^A$ acting on the $i$-th component of a tensor product).  
\begin{theo}[\cite{ESS}] There exists a unique trigonometric rational function
 $\mathcal{J}_T:\l^* \to (U_{q}(\b_-) {\otimes} U_{q}(\b_+))^\l$ such that 
\begin{enumerate}
\item $\mathcal{J}_T-q^{Z}\in I^- \otimes I^+$,
\item $\mathcal{J}_T(\lambda)$ satisfies the modified ABRR equation :
\begin{equation}\label{E:MABRR}
\mathcal{R}^{21}q_1^{2\lambda} B_1 (\mathcal{J}_T(\lambda))=\mathcal{J}_T
(\lambda)q_1^{2\lambda}q^{\Omega_\l}.
\end{equation}
\end{enumerate}
Moreover $\mathcal{J}_T(\lambda)$ satisfies the shifted 2-cocycle condition :
\begin{equation}\label{E:2cocy}
\mathcal{J}_T(\lambda)^{12,3}(\lambda)\mathcal{J}_T^{12}(\lambda+\frac{1}{2}
 h^{(3)}) =\mathcal{J}_T^{1,23}(\lambda)\mathcal{J}_T^{23}(\lambda-\frac{1}{2}
 h^{(1)}).
\end{equation}
\end{theo}
The quantum dynamical R-matrix $\widetilde{R}_T(\lambda)$ is obtained by 
twisting $\mathcal{R}$ by $\mathcal{J}_T(\lambda)$~:
$$\mathcal{R}_T(\lambda)=\mathcal{J}_T^{-1}(\lambda)\mathcal{R}^{21}
\mathcal{J}_T^{21}(\lambda),$$
$$\widetilde{R}_T(\lambda)=\mathcal{R}^{21}_T(\frac{\lambda}{\hbar}).$$
Note that the polarization we use here is the \textit{opposite} to the polarization used in \cite{ESS}, where the twist $\mathcal{J}_T(\lambda)$ was an 
element of $U_q(\b^+) \otimes U_q(\b^-)$.
\paragraph{}One aim of this paper is to provide a representation-theoretic interpretation
 of the quantum dynamical R-matrix $\mathcal{R}_T(\lambda)$. This is done in terms of twisted traces of quantum 
intertwiners and of the systems of difference equations satisfied by them.
\section{Twisted traces of quantum intertwiners}
\subsection{Definition} 
\paragraph{}Let $M_\mu$ be the Verma module over $U_q(\g)$ with highest weight
 $\mu \in \h^*$ and let $v_\mu$ be a highest weight vector. We will also
 consider the graded dual Verma module $M^*_\mu$ and let $v_\mu^*$ be its 
lowest weight vector satisfying $\langle v^*_\mu,v_\mu \rangle=1$. Let $V$ be
 a finite-dimensional $U_q(\g)$-module and let $V=\bigoplus_\nu V[\nu]$ be
 its weight space decomposition. The following result is well-known (see 
e.g \cite{ES6}) :
\begin{lem} Suppose that $M^*_\mu$ is irreducible. 
Then the map
\begin{align*}
Hom_{U_q(\g)}(M_\mu,M_\lambda \otimes V) &\to V[\mu-\lambda]\\
\Phi &\mapsto \langle v^*_\lambda,\Phi v_\mu \rangle
\end{align*}
is an isomorphism.\end{lem} 
Conversely, for every weight $\nu$ and
every homogeneous vector $v \in V[\nu]$ we will denote by
$\Phi^v_\mu:\;M_\mu \to M_{\mu-\nu}\otimes V$ the unique intertwiner
satisfying $\langle v^*_{\mu-\nu}, \Phi^v_\mu v_\mu \rangle=v$. It
will be convenient to consider all the operators $\Phi^v_\mu$
simultaneously by setting

$$\Phi^V_\mu=\sum_{v \in \mathcal{B}} \Phi^v_\mu \otimes v^* \in
Hom_\C (M_\mu, \bigoplus_\nu M_{\mu-\nu} \otimes V \otimes V^*),$$
where $\mathcal{B}$ is a homogeneous basis of $V$.

\paragraph{}Let $(\Gamma_1,\Gamma_2,T)$ be a generalized 
Belavin-Drinfeld triple. Let $\l, \h_0, C_T,\ldots$ have the same meanings as in
 Section 1. Finally, let
 $\mu,\mu'$ be weights satisfying the following relation :
\begin{equation}\label{E:mumu'}
(\mu,\alpha)=(\mu',T(\alpha))\;\mathrm{for\;all\;} \alpha \in \Gamma_1.
\end{equation}
We define a linear map $B: M_\mu \to M_{\mu'}$ by $u\cdot v_\mu \to B(u)\cdot
 v_{\mu'}$ for all $u \in U_q(\n_-)$.
\paragraph{}Now consider finite-dimensional $U_q(\g)$-modules $V_1\ldots, V_N$
 and let $v_1 \in V_1[\mu_1],$ $\ldots, v_N \in V_N[\mu_N]$ be homogeneous
 vectors such that $\overline{\mu}:=\sum_i \mu_i\in \l^\perp$. The set of pairs of weights
 $(\mu,\mu')$ satisfying (\ref{E:mumu'}) and such that $\mu'-\mu=\overline{\mu}$
 is an $\l^*$-torsor $\tilde{\l}^*$. For any such pair $(\mu,\mu')$ and for
 $\lambda \in \l^*$, we define the following formal power series in $(V_1 
\otimes
 \cdots \otimes V_N)^\l \otimes q^{2(\lambda,\mu)}\C[[q^{-2(\lambda,\alpha_1)},
 \ldots, q^{-2(\lambda,\alpha_r)}]]$ by analogy with \cite{EV4}:
$$\Psi_{v_1,\ldots, v_N}^T(\lambda,\mu)=\mathrm{Tr}_{|M_{\mu}}(\Phi^{v_1}_{\mu'
-\sum_{i=2}^N \mu_i} \cdots \Phi^{v_N}_{\mu'} B e^\lambda)$$
and
$$\Psi^T_{V_1,\ldots, V_N}(\lambda,\mu)=\sum_{v_i \in \mathcal{B}_i} 
\Psi_{v_1,\ldots, v_N}^T(\lambda,\mu)\otimes v_N^*\otimes \cdots \otimes v_1^*
,$$
where $\mathcal{B}_i$ is a homogeneous basis of $V_i$. It is clear that 
$\Psi^T_{V_1,\ldots, V_N}(\lambda,\mu)$ takes values in $(V_1 \otimes \cdots 
\otimes V_N)^\l \otimes (V_N^* \otimes \cdots \otimes V^*_1)^\l$.
\subsection{The main results}
\paragraph{}Our main result is that the functions $\Psi^T_{V_1,\ldots, V_N}(\lambda,\mu)$ satisfy some interesting difference equations. These difference equations are more conveniently expressed after some renormalizations. Set
$$J_T(\lambda)=\mathcal{J}_T(-\lambda-\rho + \frac{1}{2}(h^{(1)}+h^{(2)})),\qquad\mathbb{J}_T(\lambda)=\mathcal{J}_T(\lambda +\frac{1}{2}(h^{(1)}+h^{(2)})).$$ Put $\mathbb{Q}_T(\lambda)=m_{21}(1 \otimes S^{-1})(\mathbb{J}_T(\lambda))=m_{21}(1 \otimes S^{-1})(J_T(\lambda))$. We will denote simply by $\mathbb{Q}(\lambda)$ the element corresponding to the trivial triple $(\Gamma,\Gamma,Id)$. Also set $R_T(\lambda)=J_T^{-1}(\lambda)\mathcal{R}^{21}J_T^{21}(\lambda)$ and $\mathbb{R}_T(\lambda)=\mathbb{J}_T^{-1}(\lambda)\mathcal{R}^{21}\mathbb{J}_T^{21}(\lambda)$. Define
$$\mathbb{J}^{1\cdots N}_T(\lambda)=\mathbb{J}_T^{1,2\cdots N}(\lambda)\cdots \mathbb{J}_T^{N-1,N}(\lambda).$$
Finally, let
$$\delta_q^T(\lambda)=\big(\mathrm{Tr}_{|M_{-\rho}}(
Bq^{2\lambda})\big)^{-1}$$
be the twisted Weyl denominator. The explicit expression for $\delta_q^T(\lambda)$ is as follows. Let $\Gamma_3$ be the subset of $\Gamma_1 \cap \Gamma_2$ 
consisting of roots which return to their original position after applying 
$T$ several times, and let $\langle \Gamma_3\rangle$ be the set of positive roots which are linear combinations of roots from $\Gamma_3$. For each $\alpha \in \langle\Gamma_3\rangle$ let $N_\alpha$ be the order of the action of $B$ on $\alpha$. Consider the Lie algebra $\g$ and define $\theta_\alpha\in \C$ by $B^{N_\alpha} u_\alpha=\theta_\alpha u_\alpha$ for any $u_\alpha \in \g[\alpha]$. Then (see \cite{ES1})
$$\delta_q^T(\lambda)=q^{2(\rho,\lambda)}\prod_{\overline{\alpha} \in
\langle \Gamma_3 \rangle/B}(1-\theta_\alpha q^{-2N_\alpha(\alpha,\lambda)}).$$ Define the renormalized trace function by
$$F_{V_1,\ldots ,V_N}^T(\lambda,\mu)=[\mathbb{Q}^{-1}(\mu+h^{(*1\cdots *N)})^{(*N)}\otimes \cdots \otimes \mathbb{Q}^{-1}(\mu+h^{(*1)})^{(*1)}]\varphi^T_{V_1,\ldots,V_N}(\lambda,-\mu-\rho)$$
where
$$\varphi^T_{V_1,\ldots, V_N}(\lambda,\mu)=\mathbb{J}_T^{1\ldots, N}(\lambda)^{-1}\Psi^T_{V_1,\ldots, V_N}(\lambda,\mu)\delta_q^T(\lambda).$$
\paragraph{}Let $W$ be a finite-dimensional $U_q(\g)$-module. Consider the following difference operator acting on functions $\l^* \to (V_1\otimes \cdots \otimes V_N)^\l$~:
\begin{equation}\label{E:0pi1}
\mathcal{D}^T_W=\sum_\nu \mathrm{Tr}_{|W[\nu]}\big((\mathbb{R}_{T})^{WV_1}(\lambda + h^{(2\cdots N)})\cdots (\mathbb{R}_T)^{WV_N}(\lambda)\big)\mathbb{T}_\nu
\end{equation}
where $\mathbb{T}_\nu f(\lambda)=f(\lambda+\nu)$. In the above, we only consider the trace of the ``diagonal block'' of $(\mathbb{R}_{T})^{WV_1}(\lambda + h^{(2\cdots N)})\cdots (\mathbb{R}_T)^{WV_N}(\lambda)$, i.e the part that preserves $W[\nu]$.
\begin{theo}[Twisted Macdonald-Ruijsenaars equations]\label{T:01}
\begin{equation}\label{E:MR}
\mathcal{D}_W^T F^T_{V_1,\ldots, V_N}(\lambda,\mu)=\chi_W(q^{-2\mu})F^T_{V_1,\ldots, V_N}(\lambda,\mu),
\end{equation}
where $\chi_W(x)=\sum \mathrm{dim}\;W[\nu]x^\nu$ is the character of $W$ and where $\mathcal{D}_W^T$ acts on the variable $\lambda$.\end{theo}
\noindent
This theorem is proved in Section 3.
\paragraph{}For each $j \in \{1,\ldots,N\}$ consider the two following operators :
\begin{align}
D^T_j&=q_{*j}^{-2\mu-C_\h}q_{*j,*1}^{-2\Omega_\h}\cdots q_{*j,*j-1}^{-2\Omega_\h},\label{E:0pi21}\\
K^T_j&=\mathbb{R}_T^{j+1,j}(\lambda+h^{(j+2,\ldots, N)})^{-1} \cdots \mathbb{R}^{Nj}_T(\lambda)^{-1}\Gamma_j \mathbb{R}_T^{j1}(\lambda+h^{(2\ldots, j-1)}+h^{(j+1\ldots, N)})\times \notag\\
&\qquad \qquad \qquad \qquad \qquad \qquad \qquad \qquad \qquad \qquad  \times\cdots \mathbb{R}_T^{j,j-1}(\lambda+h^{(j+1 \ldots, N)})\label{E:0pi22}
\end{align}
where $C_\h= m_{12}(\Omega_\h)\in U(\h)$ is the quadratic Casimir element for $\h$, and where $\Gamma_jf(\lambda)=f(\lambda+h^{(j)})$.
\begin{theo}[Twisted qKZB equations]\label{T:02}
The function $F^T_{V_1,\ldots, V_N}(\lambda,\mu)$ satisfies the following difference equation for all $j=1,\ldots, N$~:
\begin{equation}\label{E:qKZB}
F^T_{V_1,\ldots, V_N}(\lambda,\mu)=(D^T_j \otimes K^T_j) F_{V_1,\ldots, V_N}^T(\lambda,\mu).
\end{equation}
\end{theo}
\noindent
This theorem is proved in Section 4.
\paragraph{}Now suppose that $(\Gamma_1,\Gamma_2,T)$ is a \textit{complete} triple, i.e $\Gamma_1=\Gamma_2=\Gamma$ and $T$ is an automorphism. In this case, the functions $F^T_V(\lambda,\mu)$ satisfy in addition some \textit{dual} difference equations, with respect to the variable $\mu$.
\paragraph{} In a complete triple the maps $B:\;U_q(\b_-) \to U_q(\b_-)$ and $B^{-1}:\;U_q(\b_+)\to U_q(\b_+)$ come from an automorphism $B: U_q(\g) \to U_q(\g)$. Let $d$ be the order of $B$ and let $\langle B \rangle \subset \mathrm{Aut}\;(U_q(\g))$ be the subgroup generated by $B$.
\paragraph{}Let $W$ be any finite-dimensional $U_q(\g)$-module. We denote by $W^B$ the twist of $W$ by $B$: as a vector space $W=W^B$ and the $U_q(\g)$-action is given by $u\cdot w=B^{-1}(g) w$. Now suppose that $W\simeq W^B$ as $U_q(\g)$-modules and let us fix an intertwiner in $\mathrm{Hom}_{U_q(\g)}(W,W^B)\subset \mathrm{Aut}_{\C}(W)$ of order $d$. This endows $W$ with the structure of a module over $\langle B \rangle \ltimes U_q(\g)$.
\paragraph{}Consider the following difference operator acting on functions with values in $(V_N^* \otimes\cdots \otimes V^*_1)^\l$ :
\begin{equation}\label{E:70}
 \mathcal{D}_W^{\vee,T}=\sum_{\nu} \mathrm{Tr}_{|W[\nu]}\big(\mathbb{R}^{WV^{*}_N}(\mu + h^{(*1\cdots *N-1)})\cdots \mathbb{R}^{WV^*_1}(\mu)B_{W}\big)\mathbb{T}^\vee_\nu
\end{equation}
 where $\mathbb{T}^{\vee}_\nu f(\mu)=f(\mu+\nu)$.
\begin{theo}[Dual twisted Macdonald-Ruijsenaars equations]\label{T:03} 
\begin{equation}\label{E:DMR}
\mathcal{D}^{\vee,T}_{W}F^T_{V_1,\ldots, V_N}(\lambda,\mu)=\mathrm{Tr}_{|W^{\h_0}}(q^{-2\lambda}B)F^T_{V_1,\ldots, V_N}(\lambda,\mu).
\end{equation}
\end{theo}
\paragraph{Remark.} Any $B$-invariant finite-dimensional $U_q(\g)$ module is a direct sum of modules $\overline{V}_{\nu_0}:=\bigoplus_{\nu \in \langle B \rangle \nu_0} V_\nu$ where $V_\nu$ is the irreducible highest weight module of highest weight $\nu$ and where $\nu_0$ is dominant integral. It is easy to see that both sides of (\ref{E:DMR}) identically vanish when $W=\overline{V}_{\nu_0}$ and $\nu_0 \not\in \l^*$ (i.e when $B(\nu_0) \neq \nu_0$). 
\paragraph{}The twisted character $\mathrm{Tr}_{|W^{\h_0}}(q^{2\lambda}B)$ can be expressed explicitly when $W=V_\nu$ with $\nu \in \l^*$. Consider the quotient Dynkin diagram $\overline{\Gamma}=\{\overline{\alpha}=\sum_{\alpha \in \overline{\alpha}} \alpha\,|\, \overline{\alpha} \in \Gamma/\langle T\rangle\}$ which forms a set of simple roots in $\l^*$ with respect to the restriction of $(\,,\,)$ to $\l^*$. Namely, if $\Gamma=\{\alpha_1,\ldots, \alpha_{n}\}$ is of type $A_{n}$ and if $T$ is the ``flip'' $T: \alpha_i \mapsto \alpha_{n+1-i}$ then $\Gamma/\langle T \rangle$ is $B_{k}$ where $n=2k-1$ or $n=2k$; if $\Gamma=D_4$ and $T$ is the rotation of order 
three around the trivalent root then $\Gamma/\langle T\rangle=G_2$; if $\Gamma=D_k$ and $T$ is the symmetry of order two around the trivalent root then $\overline{\Gamma}=C_{k-1}$; and if $\Gamma=E_6$ and $T$ is the symmetry around the trivalent root then $\Gamma/\langle T \rangle=F_4$. Let $\overline{\Delta}$ be the root system of $\overline{\Gamma}$. Let $\overline{\g}$ be the simple complex Lie algebra associated to $\overline{\Gamma}$. Note that any weight $\nu \in \l^*$ is naturally a weight for $\overline{\Gamma}$. However, the scalar product on $\l^*$ is \textit{not} the usual one corresponding to the root system $\overline{\Delta}$. For instance we have $(\overline{\alpha},\overline{\alpha})=2N_\alpha$ if $N_\alpha$ is the number of elements in the $T$-orbit $\overline{\alpha}$ and if any two elements of that orbit are orthogonal. For $\lambda=\sum_\alpha c_{\overline{\alpha}} \overline{\alpha} \in \l^*$ set $\overline{\lambda}=\sum_{\overline{\alpha}} \frac{2N_\alpha}{(\overline{\alpha},\overline{\alpha})}c_{\overline{\alpha}} \overline{\alpha}$.
\begin{prop}\label{P:001} For any $\nu \in \l^*$ we have
$$\mathrm{Tr}_{|W^{\h_0}}(q^{2\lambda}B)=\chi_{\overline{V}_{\nu}}(q^{2\overline{\lambda}})$$
where $\overline{V}_\nu$ is the irreducible $\overline{\g}$-module of highest weight $\nu$.\end{prop}
\noindent
Theorem~\ref{T:03} and Proposition~\ref{P:001} are proved in Section 5.
\paragraph{}Now, define for each $j \in \{1,\ldots,N\}$ the following operators
\begin{equation}\label{E:71}
\begin{split}
D^{\vee,T}_j&=q_j^{-2\lambda-C_\l}q_{j,j+1}^{-\Omega_\l}\cdots q_{j,N}^{-\Omega_\l},\\
K^{\vee,T}_j&=\mathbb{R}_{*j-1,*j}(\mu+h^{(*1\cdots *j-2)})^{-1}\cdots \mathbb{R}_{*1,*j}(\mu)^{-1}\Gamma^*_{B^{-1}(j)}\times\\
\mathbb{R}_{*j,*N}(\mu&+h^{(*j+1\cdots *N-1)}+h^{(*1\cdots *j-1)})
\cdots \mathbb{R}_{*j,*j+1}(\mu+h^{(*1\cdots *j-1)}),
\end{split}
\end{equation}
where $C_\l=m_{12}(\Omega_\l) \in U(\l)$ and where $\Gamma^*_{B^{-1}(j)}f(\mu)=f(\mu+B^{-1}(h^{(*j)}))$.
\begin{theo}[Dual twisted qKZB equations]\label{T:04}
The functions $F^T_{V_1,\ldots, V_N}(\lambda,\mu)$ satisfy the following difference equation for each $j=1\ldots, N$ :
\begin{equation}\label{E:DqKZB}
B_{V_j}B^*_{V_j^*}F_{V_1,\ldots, V_j,\ldots, V_N}^T(\lambda,\mu)=(D^{\vee,T}_j\otimes K^{\vee,T}_j) F^T_{V_1,\ldots, V_j^B,\ldots, V_N}(\lambda,\mu).
\end{equation}
\end{theo}
\noindent
This theorem is proved in Section 6.  
\paragraph{Remark 1.} For $T=Id$, Theorems~\ref{T:01}-\ref{T:04} appear in \cite{EV4}.
\paragraph{Remark 2.} We do not expect the dual equations to exist for non-complete triples. This can be explained in the following way. Suppose that $\g$ is an affine Lie algebra and that $T=Id$, so that $r_T(\lambda,z)$ is the Felder elliptic dynamical r-matrix. In that case it is known that the the \textit{dual} trigonometric qKZB equations \textit{without spectral parameter} can be interpreted as monodromy of the flat connection on the torus defined by the classical (elliptic) KZB equations (see \cite{Kir}). One can show that this is true for any elliptic dynamical r-matrix. On the other hand, it was proved in \cite{ES1} Proposition 4.2 that the classical dynamical r-matrix with spectral parameter $r_T(\lambda,z)$ associated to an affine Lie algebra and a triple $(\Gamma_1,\Gamma_2,T)$ is \textit{elliptic} only when $T$ is an automorphism; for general triples, it is partially elliptic and partially trigonometric (for instance, it is purely trigonometric when $T$ is nilpotent). This shows that the monodromy of these KZB equations should be defined only for complete triples, and hence the existence of the dual equations should be expected only for them.
\paragraph{Remark 3.} The above theorems are also valid for the specialized quantum group $U_q(\g)$, which is obtained from the formal quantum group when we take $\hbar\in \C^*\setminus\{i\mathbb{R}\}$ to be a complex number. In that case, it is more convenient to consider the twist $\mathcal{J}_T(\lambda)$ as an endomorphism of the functor
$$F:\;Rep(U_q(\g)) \times Rep(U_q(\g)) \to Vec$$
which assigns to any two finite-dimensional $U_q(\g)$-modules $V$ and $W$ the vector space $V \otimes W$. Here $Rep(U_q(\g))$ is the category of finite-dimensional $U_q(\g)$-modules and $Vec$ is the category of finite-dimensional $\C$-vector spaces. For instance, equation (\ref{E:2cocy}) means that for every three representations $U,V,W$ and vectors $u \in U$, $v \in V$ and $w \in W$ with respective weights $\lambda_u,\lambda_v$ and $\lambda_w$ we have
$$ \mathcal{J}_T(\lambda)^{12,3}(\lambda)\mathcal{J}_T^{12}(\lambda+\frac{1}{2}
 \lambda_w)(u \otimes v \otimes w) =\mathcal{J}_T^{1,23}(\lambda)\mathcal{J}_T^{23}(\lambda-\frac{1}{2}
\lambda_u)(u \otimes v \otimes w).$$
\section{The twisted Macdonald-Ruijsenaars equations}
\paragraph{}The proof of Theorem~\ref{T:01} is an extension of the proof of Theorem 1.1
 of \cite{EV4} to the case of an arbitrary generalized 
Belavin-Drinfeld triple. From now on we fix such a triple $(\Gamma_1,
\Gamma_2,T)$. We first note that the notion of radial part 
generalizes straightforwardly to the twisted setting :
\begin{prop} Let $V$ be a finite-dimensional $U_q(\g)$-module. For any $X \in
 U_q(\g)$ there exists a unique difference operator
 $\mathcal{D}^T_X$ (with respect to the variable $\lambda$) acting on formal 
power series in $V^\l \otimes q^{2(\lambda,\mu)} \C[[q^{-2(\lambda,\alpha_1)}
,\ldots, ,q^{-2(\lambda,\alpha_r)}]]$, $\lambda \in \l^*$ such that we have  
$$\mathrm{Tr}_{|M_\mu}(\Phi^V_{\mu'} X q^{2\lambda}B)=
\mathcal{D}_X^T \mathrm{Tr}_{|M_\mu}(\Phi^V_{\mu'} q^{2\lambda}B).$$
\end{prop}
The operator $\mathcal{D}_X^T$ is called the twisted radial part of $X$.
\paragraph{}For any finite-dimensional $U_q(\g)$-module $W$ set
$$C_W=\mathrm{Tr}_{|W}(1 \otimes \pi_W)(\mathcal{R}^{21}\mathcal{R}(1
 \otimes q^{2\rho})).$$
It is well-known (see \cite{D}, \cite{R}) that the map $W \to C_W$ defines a
 homomorphism from the Grothendieck ring of the category of finite-dimensional
 $U_q(\g)$-modules to the center of $U_q(\g)$. Set
 $\mathcal{M}^T_W=\mathcal{D}^T_{C_W}$.
\begin{prop}\label{P:02} We have
\begin{enumerate}
\item $\mathcal{M}_W^T \mathcal{M}_V^T=\mathcal{M}_V^T \mathcal{M}_W^T$ for
 all $V,W$,
\item $\mathcal{M}^T_W \Psi^T_V(\lambda,\mu)=\chi_W(q^{2(\mu +\rho)}) 
\Psi^T_V(\lambda,\mu)$ where $\chi_W(x)=\sum_\nu \mathrm{dim}\;W[\nu] x^\nu$
 is the character of $W$.
\end{enumerate}
\end{prop}
\noindent
\textit{Proof.} See \cite{EKi}, \cite{EV4}.
\paragraph{}Let us now proceed to explicitly compute the operator
 $\mathcal{M}_W$.
\paragraph{}Let $V$ be a finite-dimensional $U_q(\g)$-module. Introduce the
 following function with values in $V \otimes V^* \otimes U_q(\g)$, with 
components labeled as $1, *1$ and $2$ respectively~: 
$$Z_V(\lambda,\mu)=\mathrm{Tr}_{|M_\mu}(\Phi_{\mu'}^{V} \mathcal{R}^{20}B_0 
q_0^{2\lambda}).$$
\begin{lem}\label{L:01} We have
\begin{equation}\label{E:01}
Z_V(\lambda,\mu)=\mathcal{J}_T^{12}(\lambda)\Psi^T_V(\lambda+\frac{1}{2}
h^{(2)},\mu).
\end{equation}
\end{lem}
\noindent
\textit{Proof.} First we note that, by pulling the R-matrix around the trace
 and using the intertwining property together with the fact that $B_1^{-1}\mathcal{R}=B_2\mathcal{R}$ we obtain
\begin{equation*}
\begin{split}
Z_V(\lambda,\mu)&=\mathcal{R}^{21}q_1^{2\lambda}\mathrm{Tr}_{|M_\mu}
(\Phi_{\mu'}^V (B_2^{-1}\mathcal{R}^{20}) B_0 q_0^{2\lambda})\\
&=\mathcal{R}^{21}q_1^{2\lambda} B_2^{-1}\big(\mathrm{Tr}_{|M_\mu}(\Phi_{\mu'}^V
\mathcal{R}^{20} B_0 q_0^{2\lambda})\big)\\
&=\mathcal{R}^{21}q_1^{2\lambda} B_2^{-1}Z_V(\lambda,\mu)
\end{split}
\end{equation*}
 On the other hand, using the defining equation for $\mathcal{J}_T(\lambda)$, the relation $B_2^{-1}\mathcal{J}_T(\lambda)=B_1
\mathcal{J}_T(\lambda)$ and the $\l$-invariance of $\Psi_V^T(\lambda,\mu)$ we have
\begin{equation*}
\begin{split}
\mathcal{R}^{21}q_1^{2\lambda}B_2^{-1}\big(\mathcal{J}_T^{12}(\lambda)\Psi^T_V
(\lambda+\frac{1}{2}h^{(2)},\mu)\big)
=&\mathcal{R}^{21}q_1^{2\lambda}\big(B_1\mathcal{J}_T^{12}(\lambda)\big)
\Psi^T_V(\lambda+\frac{1}{2}h^{(2)},\mu)\\
=&\mathcal{J}_T^{12}(\lambda)q_1^{2\lambda}q_{12}^{\Omega_\l}\Psi^T_V(\lambda+
\frac{1}{2}h^{(2)},\mu)\\
=&\mathcal{J}_T^{12}(\lambda)\Psi^T_V(\lambda+\frac{1}{2}h^{(2)},\mu)
\end{split}
\end{equation*}
 The lemma now follows from the 
fact both sides of (\ref{E:01}) satisfy the relation $Y= \mathcal{R}^{21}
q_1^{2\lambda} B_2^{-1}Y$ and are of the form $Y=q^Z\Psi^T_V(\lambda +
\frac{1}{2}h^{(2)},\mu) + l.o.t$, where $l.o.t$ stands for terms of strictly positive degree in component $2$.\qed
\paragraph{}Consider the following function with values in $V \otimes V^* 
\otimes U_q(\g) \otimes U_q(\g)$ (with components labeled as $1,*1,2$ and $3$
 respectively) :
$$X_V(\lambda,\mu)=\mathrm{Tr}_{|M_\mu}(\Phi_{\mu'}^V \mathcal{R}^{20}q_0^{2
\lambda}B_0(\mathcal{R}^{03})^{-1}).$$
\begin{lem}We have
\begin{equation}\label{E:02}
X_V(\lambda,\mu)=\mathcal{J}_T^{3,12}(\lambda)\mathcal{J}_T^{12}(\lambda-
\frac{1}{2}h^{(3)})\Psi^T_V(\lambda+\frac{1}{2}(h^{(2)}-h^{(3)}),\mu)
\mathcal{J}_T^{32}(\lambda)^{-1}
\end{equation}
\end{lem}
\noindent
\textit{Proof.}
Moving $\mathcal{R}^{03}$ around the trace, using the quantum Yang-Baxter equation for $\mathcal{R}$ and the $B$-invariance property of $\mathcal{R}$ again, we get
\begin{equation*}
\begin{split}
X_V(\lambda,\mu)&=\mathcal{R}^{13}\mathrm{Tr}_{|M_\mu}(\Phi_{\mu'}^V(
\mathcal{R}^{03})^{-1}\mathcal{R}^{20}q_0^{2\lambda}B_0)\\
&=\mathcal{R}^{13}\mathcal{R}^{23}\mathrm{Tr}_{|M_\mu}(\Phi_{\mu'}^V
\mathcal{R}^{20}(\mathcal{R}^{03})^{-1}q_0^{2\lambda}B_0)
(\mathcal{R}^{23})^{-1}\\
&=\mathcal{R}^{13}\mathcal{R}^{23}q_3^{2\lambda}\mathrm{Tr}_{|M_\mu}
(\Phi_{\mu'}^V
\mathcal{R}^{20}q_0^{2\lambda}(\mathcal{R}^{03})^{-1}B_0)q_3^{-2\lambda}
(\mathcal{R}^{23})^{-1}\\
&=\mathcal{R}^{13}\mathcal{R}^{23}q_3^{2\lambda}B_3\big(\mathrm{Tr}_{|M_\mu}
(\Phi_{\mu'}^V\mathcal{R}^{20}q_0^{2\lambda}B_0(\mathcal{R}^{03})^{-1})\big)
q_3^{-2\lambda}(\mathcal{R}^{23})^{-1}\\
&=\mathcal{R}^{13}\mathcal{R}^{23}q_3^{2\lambda}B_3\big(X_V(\lambda,\mu)\big)
q_3^{-2\lambda}(\mathcal{R}^{23})^{-1}.
\end{split}
\end{equation*}
 On the other hand, using the modified ABRR equation (\ref{E:MABRR}) we have
\begin{equation*}
\begin{split}
\mathcal{R}^{13}&\mathcal{R}^{23}q_3^{2\lambda}B_3\big(\mathcal{J}_T^{3,12}
(\lambda)\big)\mathcal{J}_T^{12}(\lambda-\frac{1}{2}h^{(3)})\Psi^T_V(\lambda+
\frac{1}{2}(h^{(2)}-h^{(3)}),\mu)\times\\
&\qquad \qquad \qquad \qquad \qquad \qquad\qquad \qquad\qquad 
\times B_3(\mathcal{J}_T^{32}(\lambda)^{-1})q_3^{-2\lambda}
(\mathcal{R}^{23})^{-1}\\
&=\Delta_1\big(\mathcal{R}^{13}q_3^{2\lambda}B_3(\mathcal{J}_T^{31}(\lambda))
\big)\mathcal{J}_T^{12}(\lambda-\frac{1}{2}h^{(3)})\Psi^T_V(\lambda+
\frac{1}{2}(h^{(2)}-h^{(3)}),\mu)\times
\\&\qquad \qquad \qquad \qquad\qquad \qquad \qquad \qquad\qquad \qquad\qquad 
\times q_3^{-2\lambda}q_{23}^{-\Omega_\l}\mathcal{J}_T^{32}(\lambda)^{-1}\\
&=\Delta_1\big(\mathcal{J}_T^{31}(\lambda)q_3^{2\lambda}q_{31}^{\Omega_\l}
\big)\mathcal{J}_T^{12}(\lambda-\frac{1}{2}h^{(3)})\Psi^T_V(\lambda+
\frac{1}{2}(h^{(2)}-h^{(3)}),\mu)\times\\
&\qquad \qquad \qquad \qquad\qquad \qquad \qquad \qquad\qquad \qquad\qquad
 \times q_3^{-2\lambda}q_{23}^{-\Omega_\l}\mathcal{J}_T^{32}(\lambda)^{-1}\\
&=\mathcal{J}_T^{3,12}q_3^{2\lambda}q_{31}^{\Omega_\l}q_{32}^{\Omega_\l}
\mathcal{J}_T^{12}(\lambda-\frac{1}{2}h^{(3)})\Psi^T_V(\lambda+\frac{1}{2}
(h^{(2)}-h^{(3)}),\mu)q_3^{-2\lambda}q_{23}^{-\Omega_\l}\mathcal{J}_T^{32}
(\lambda)^{-1}\\
&=\mathcal{J}_T^{3,12}\mathcal{J}_T^{12}(\lambda-\frac{1}{2}h^{(3)})\Psi^T_V
(\lambda+\frac{1}{2}(h^{(2)}-h^{(3)}),\mu)\mathcal{J}_T^{32}(\lambda)^{-1}
\end{split}
\end{equation*}
Now set $X(\lambda)=
(\mathcal{J}_T^{3,12})^{-1}X_V(\lambda)\mathcal{J}^{32}_T(\lambda)$. By the 
above and by Lemma~\ref{L:01}, both $X(\lambda)$ and $Z_V(\lambda-\frac{1}{2}
h^{(3)})$ satisfy the equation
$$q_3^{2\lambda}q_{31}^{\Omega_\l}q_{32}^{\Omega_\l}Y=Yq_3^{2\lambda}
q_{23}^{\Omega_\l}$$
and are both of the form $Y=\mathcal{J}_T^{12}(\lambda-\frac{1}{2}h^{(3)})
\Psi_V^T(\lambda+\frac{1}{2}(h^{(2)}-h^{(3)}))+l.o.t$. Hence $X(\lambda)= Z_V
(\lambda-\frac{1}{2}h^{(3)})$ and the lemma is proved.\qed
\begin{cor}\label{C:01} We have
\begin{equation*}
\begin{split}
\mathrm{Tr}_{|M_\mu}&(\Phi_\mu^V\mathcal{R}^{20}(\mathcal{R}^{03})^{-1}
q^{2\lambda}B_0)\\
=&B_3\bigg(q_3^{2\lambda}\mathcal{J}_T^{3,12}(\lambda)\mathcal{J}_T^{12}
(\lambda-\frac{1}{2}h^{(3)})\Psi_V(\lambda+\frac{1}{2}(h^{(2)}-h^{(3)}),\mu)
\mathcal{J}^{32}_T(\lambda)^{-1}q_3^{-2\lambda}\bigg)
\end{split}
\end{equation*}
\end{cor}
\paragraph{}Now let $W$ be any finite-dimensional $U_q(\g)$-module. By 
Corollary~\ref{C:01}, we get
\begin{equation}\label{E:025}
\begin{split}
\mathcal{M}_W^T\Psi_V^T(\lambda,\mu)&=\mathrm{Tr}_{|M_\mu}(\Phi_{\mu'}^VC_W
q_0^{2\lambda}B_0)\\
&=\mathrm{Tr}_{|M_\mu}\big(\mathrm{Tr}_{|W}(\mathcal{R}^{W0}\mathcal{R}^{0W}
q_W^{2\rho})\Phi_{\mu'}^Vq_0^{2\lambda}B_0)\\
&=\mathrm{Tr}_{|W}\mathrm{Tr}_{|M_\mu}\big(m_{23}(\mathcal{R}^{20}
S_3(\mathcal{R}^{03})^{-1}q_3^{2\rho})\Phi_{\mu'}^Vq_0^{2\lambda}B_0\big)\\
&=\mathrm{Tr}_{|W}\big\{m_{23}\big(S_3\mathrm{Tr}_{|M_\mu}(\mathcal{R}^{20}
(\mathcal{R}^{03})^{-1}\Phi_\mu^Vq_0^{2\lambda}B_0)\big)q_2^{2\rho}\big\}\\
&=\mathrm{Tr}_{|W}\big\{\sum_{ijk}[d_k^{(1)}(\lambda)c_i\otimes q^{-2\lambda}
d_k^{(2)}(\lambda)d_i(\lambda+\frac{1}{2}h^{(2)})]\times\\
&\qquad\qquad\times\Psi_V^T(\lambda+h^{(W)}) \big(d'_j
(\lambda)q^{2\lambda}S(B(c'_j))S(B(c_k))q^{2\rho}\big)_W\big\}
\end{split}
\end{equation}
where $m: U_q(\g) \otimes U_q(\g)\to U_q(\g)$ is the multiplication map,
$\mathcal{J}_T(\lambda)=\sum_i c_i \otimes d_i(\lambda)$, $\mathcal{J}_T^{-1}
(\lambda)=\sum_ic'_i \otimes d'_i(\lambda)$ and where we used Sweedler's 
notation for coproducts : $\Delta(x)=\sum x^{(1)} \otimes x^{(2)}$. 
\paragraph{}Let us set $\mathcal{R}=\sum_i a_i \otimes b_i$, $\mathcal{R}^{-1}=\sum_i a'_i 
\otimes b'_i$. 
\begin{lem} We have
\begin{equation}\label{E:03}
\sum_j d_j'(\lambda)q^{2\lambda}S(B(c'_j))=q^{C_\l}P(\lambda)
S(u)q^{2\lambda}
\end{equation}
where $u=\sum_i S(b_i)a_i$ is the Drinfeld element and $P(\lambda)=
\sum_j d'_j(\lambda)S^{-1}(c'_j)$.
\end{lem}
\noindent
\textit{Proof.} This is obtained by applying $m_{21}(S \otimes 1)$ to the relation $(B_1\mathcal{J}_T
(\lambda))^{-1}q_1^{-2\lambda}=q^{-\Omega_\l}q_1^{-2\lambda}\mathcal{J}_T^{-1}
(\lambda)\mathcal{R}^{21}$, which itself follows from (\ref{E:MABRR}) and the
 $B$-invariance of $\mathcal{J}_T(\lambda)$.\qed
\paragraph{}Substituting (\ref{E:03}) in (\ref{E:025}) yields
\begin{equation}\label{E:04}
\begin{split}
\mathcal{M}_W^T\Psi_V^T(\lambda,\mu)&=\sum_{ijk,\nu} d_k^{(1)}(\lambda)c_i
\mathrm{Tr}_{|W[\nu]}\big\{q^{C_\l}P(\lambda)q^{2\lambda}S^{-1}
(B(c_k))q^{-2\lambda} d_k^{(2)}(\lambda)\times\\
&\qquad \qquad \qquad \qquad \times d_i(\lambda+\frac{1}{2}\nu)S(u)
q^{2\rho}\big\}\Psi_V^T(\lambda+\nu).
\end{split}
\end{equation}
\begin{lem}~\label{L:06} We have
\begin{equation}\label{E:05}
\begin{split}
\sum_k &d_k^{(1)}(\lambda)\otimes q^{2\lambda}S^{-1}(B(c_k))q^{-2\lambda}
d_k^{(2)}(\lambda)\\
&= \sum_{jk} (a'_j)^{(1)}d_k^{(1)}(\lambda) q^{-\Omega_\l+1\otimes C_\l}
\big\{S^{-1}(c_k)S^{-1}(b'_j)(a'_j)^{(2)}d_k^{(2)}(\lambda)\}_2
\end{split}
\end{equation}
\end{lem}
\noindent
\textit{Proof.} From the modified ABRR relation we get
$$q_1^{2\lambda}B_1(\mathcal{J}_T(\lambda))q_1^{-2\lambda}=
(\mathcal{R}^{21})^{-1}\mathcal{J}_T(\lambda)q^{\Omega_\l}.$$
Applying $1 \otimes \Delta$ yields
$$q_1^{2\lambda}B_1(\mathcal{J}_T^{1,23}(\lambda))q_1^{-2\lambda}=
(\mathcal{R}^{23,1})^{-1}\mathcal{J}_T^{1,23}(\lambda)q^{\Omega_\l}_{12}
q^{\Omega_\l}_{13},$$
which can be written as
\begin{equation*}
\begin{split}
\sum_k q_1^{2\lambda}B(c_k)q_1^{-2\lambda}& \otimes d_k^{(1)}(\lambda)
 \otimes d_k^{(2)}(\lambda)\\
&=\sum_{ik} (b'_i \otimes (a'_i)^{(1)} \otimes (a'_i)^{(2)})\times (c_k
 \otimes d_k^{(1)}(\lambda)\otimes d_k^{(2)}(\lambda)) q_{12}^{\Omega_\l}
q_{13}^{\Omega_\l}.
\end{split}
\end{equation*}
Equation (\ref{E:05}) is now obtained by applying $m_{13}(S^{-1} \otimes 1 
\otimes 1)$.  \qed
\paragraph{}We introduce the following notation. For any linear operator $H(\lambda)\in\\ \mathrm{End}\;(V_1 \otimes \cdots \otimes V_N)$ we set
$$H(\lambda+\hat{h}^{(i)})(v_1\otimes \cdots \otimes v_N)=\sum_{\nu}H_\nu(\lambda+\nu)(v_1 \otimes \cdots \otimes v_N)$$
where $H_\nu(\lambda):V_1 \otimes \cdots \otimes V_i \otimes \cdots \otimes V_N \to V_1 \otimes \cdots \otimes V_i[\nu] \otimes \cdots \otimes V_N$ is the block of $H(\lambda)$ with image $V_i[\nu]$ in the $i$-th component. In other words, we replace $\hat{h}^{(i)}$ by the weight in the $i$-th component \textit{after} the action of $H$.
\begin{lem} The following identities hold :
\begin{enumerate}
\item[i)] $\sum_j (a'_j)^{(1)} \otimes S^{-1}(b'_j)(a'_j)^{(2)}=\sum a_k
 \otimes u b_k,$
\item[ii)] $\mathcal{R}^{23}\mathcal{J}_T^{1,23}(\lambda)=\mathcal{J}_T^{1,32}
(\lambda)\mathcal{R}^{23},$
\item[iii)] $\sum S(c_i)d_i^{(1)}(\lambda)\otimes d_i^{(2)}(\lambda)=S(\mathbb{Q}_T)
(\lambda-\frac{1}{2}h^{(2)})_1\mathcal{J}_T^{-1}(\lambda +\frac{1}{2}
\hat{h}^{(1)}).$
\end{enumerate}
\end{lem}
\noindent
\textit{Proof.} Equalities i) and iii) are proved in the same fashion as in \cite{EV4}. Equality ii) follows from the relation $\mathcal{R}\Delta=\Delta^{op}\mathcal{R}$.\qed
%To
% prove ii), note that by the modified ABRR relation,
%$$B_1(\mathcal{J}_T^{1,23})(\lambda)=q_1^{-2\lambda}(\mathcal{R}^{31})^{-1}
%(\mathcal{R}^{21})^{-1}\mathcal{J}_T^{1,23}(\lambda)q_1^{2\lambda}
%q_{12}^{\Omega_\l}q_{13}^{\Omega_\l}$$
%hence
%$$B_1(\mathcal{J}_T^{1,32}(\lambda))=q_1^{-2\lambda}(\mathcal{R}^{21})^{-1}
%(\mathcal{R}^{31})^{-1}\mathcal{J}_T^{1,32}(\lambda)q_1^{2\lambda}
%q_{12}^{\Omega_\l}q_{13}^{\Omega_\l}$$
%which implies that both $\mathcal{J}_T^{1,23}(\lambda)$ and
% $(\mathcal{R}^{23})^{-1}\mathcal{J}_T^{1,32}(\lambda)\mathcal{R}^{23}$ 
%satisfy the equation 
%$$B_1Y=q_1^{-2\lambda}(\mathcal{R}^{31})^{-1}(\mathcal{R}^{21})^{-1}Y
%q_1^{2\lambda}q_{12}^{\Omega_\l}q_{13}^{\Omega_\l}$$
%and are of the form $Y=q^{Z_{12}+Z_{13}}+l.o.t$. Thus ii) holds.
\begin{cor}\label{C:02} We have
\begin{equation*}
\begin{split}
\sum d_k^{(1)}&(\lambda)\otimes q^{2\lambda}S^{-1}(B(c_k))q^{-2\lambda}
d_k^{(2)}(\lambda)\\
&=q^{-\Omega_\l -1\otimes C_\l}u_2^{-1}S(\mathbb{Q}_T)_2(\lambda-\frac{1}{2}
h^{(1)})(\mathcal{J}_T^{21})^{-1}(\lambda+\frac{1}{2}\hat{h}^{(2)})\mathcal{R}.
\end{split}
\end{equation*}
\end{cor}
\noindent
\textit{Proof.} Use i),ii) and iii) successively, as in \cite{EV4}, (2.32).\qed
\paragraph{}By Corollary~\ref{C:02} and using the relation $B_1\mathcal{J}_T(\lambda)=B_2^{-1}\mathcal{J}_T(\lambda)$, we can rewrite (\ref{E:04}) as follows :
\begin{equation}\label{E:06}
\begin{split}
&\mathcal{M}_W^T\Psi_V^T(\lambda,\mu)\\
&=\sum_\nu \mathrm{Tr}_{|W[\nu]}\bigg\{P_2(\lambda)
q^{-\Omega_\l-1\otimes C_\l}u_2^{-1}S(\mathbb{Q}_T)_2(\lambda-\frac{1}{2}
h^{(1)})(\mathcal{J}_T^{21})^{-1}(\lambda+\frac{1}{2}\hat{h}^{(2)})
\mathcal{R}\\
&\qquad \qquad \qquad \qquad \qquad \quad \times (c_i)_1d_i(\lambda+
\frac{1}{2}\nu)_2 S(u)_2q_2^{2\rho}q_2^{m_{12}\Omega_\l}\bigg\}\Psi_V^T(\lambda+
\nu,\mu)\\
&=\sum_\nu\mathrm{Tr}_{|W[\nu]}\big(\tilde{G}(\lambda)(\mathbb{R}_T)^{WV}
(\lambda)\big)\Psi_V(\lambda+\nu,\mu)
\end{split}
\end{equation}
where $\tilde{G}(\lambda)=q^{-2\rho}P(\lambda)S(\mathbb{Q}_T)(\lambda)$. We now proceed to compute $\tilde{G}(\lambda)$. 
\begin{prop}\label{P:03}We have $\tilde{G}(\lambda)=\frac{\delta^T_q
(\lambda+h)}{\delta_q^T(\lambda)}$.
\end{prop}
\noindent
\textit{Proof.} The following lemma is proved as in \cite{EV4}~:
\begin{lem}We have $P(\lambda)=\mathbb{Q}_T^{-1}(\lambda+h)$, i.e $\tilde{G}(\lambda)=
G(\lambda+h)$ where $G(\lambda)=q^{-2\rho}\mathbb{Q}_T^{-1}(\lambda)S(\mathbb{Q}_T)(\lambda-h)$.
\end{lem}
\noindent
A direct (though lengthy) computation shows that
\begin{equation}\label{E:07}
\Delta(G(\lambda))=\mathbb{J}_T(\lambda)\big(G(\lambda+h^{(2)})\otimes G(\lambda)\big)
 \mathbb{J}_T^{-1}(\lambda)
\end{equation}
(see \cite{Majid} for a detailed proof of this in the nondynamical case; the dynamical case is analogous). In
 particular, replacing $\lambda$ by $\frac{\lambda}{\hbar}$, we have
\begin{equation*}
\begin{split}
\mathcal{R}\mathbb{J}_T(\frac{\lambda}{\hbar})(G(\frac{\lambda+
\hbar h^{(2)}}{\hbar})\otimes G(\frac{\lambda}{\hbar}))&\mathbb{J}_T^{-1}
(\frac{\lambda}{\hbar})\\
&=\mathbb{J}_T^{21}(\frac{\lambda}{\hbar})(G(\frac{\lambda}{\hbar})\otimes
 G(\frac{\lambda+\hbar h^{(1)}}{\hbar}))(\mathbb{J}^{21}_T)^{-1}
(\frac{\lambda}{\hbar})\mathcal{R}
\end{split}
\end{equation*}
which can be rewritten as
\begin{equation}\label{E:08}
\mathbb{R}^{21}_T(\frac{\lambda}{\hbar})(G(\frac{\lambda+
\hbar h^{(2)}}{\hbar})\otimes G(\frac{\lambda}{\hbar}))=
(G(\frac{\lambda}{\hbar})\otimes G(\frac{\lambda+\hbar h^{(1)}}{\hbar}))
\mathbb{R}^{21}_T(\frac{\lambda}{\hbar})
\end{equation}
Let us now expand $\mathbb{R}_T(\frac{\lambda}{\hbar})$ and 
$G(\frac{\lambda}{\hbar})$ around $\hbar=0$~:
$$\mathbb{R}_T(\frac{\lambda}{\hbar})=1+\hbar r(\lambda) +
\mathcal{O}(\hbar^2),\qquad G(\frac{\lambda}{\hbar})=1+\hbar g_1(\lambda) +
 \mathcal{O}(\hbar^2)$$
where $g_1(\lambda)=(G(\lambda/\hbar)-1)/\hbar \in U_q(\g)/\hbar U_q(\g) \simeq U(\g)$.
Note that by (\ref{E:07}) we have $\Delta_0(g_1(\lambda))=g_1(\lambda) \otimes
 1 + 1 \otimes g_1(\lambda)$ (where $\Delta_0$ is the usual coproduct on $U\g$), which implies that $g_1(\lambda) \in \g$. But
 since $G(\lambda)$ is of $\l$-weight zero, $g_1(\lambda) \in \h$. Now, by
 (\ref{E:08}), we have
$$\sum_i x_i \wedge \frac{\partial g_1(\lambda)}{\partial x_i}=[r(\lambda),
 g_1(\lambda) \otimes 1 + 1 \otimes g_1(\lambda)],$$
where $(x_i)$ is a basis of $\l$. In particular, $[r(\lambda), g_1(\lambda)
 \otimes 1 + 1 \otimes g_1(\lambda)] \in \Lambda^2\h$. But this implies that
 $[r(\lambda), g_1(\lambda) \otimes 1 + 1 \otimes g_1(\lambda)]=0$. Thus $g_1:\l^* \to \l$ is a closed 1-form on $\l^*$ and there exists functions $f_1(\lambda)$
 and $g_2(\lambda)=\frac{1}{\hbar^2}(G(\frac{\lambda}{\hbar})-\frac{f_1(\frac{\lambda}{\hbar})}{f_1
(\frac{\lambda-\hbar h}{\hbar})}) \in U_q(\g)/\hbar U_q(\g)\simeq U(\g)$.
By the same argument, $\g_2(\lambda)$ is a closed 1-form. Continuing in this
 process, we finally obtain a function $f$ defined on $\l^*$ such that 
$G(\lambda)=
\frac{f(\lambda)}{f(\lambda-h)}$. It remains to determine
 $f(\lambda)$ explicitly. For this, apply (\ref{E:06}) and Proposition~\ref{P:02} 2. to
 the case of the trivial representation $V=\C$. Then $\Psi^T_V(\lambda,\mu)=
\frac{q^{2(\mu+\rho,\lambda)}}{\delta_q^T(\lambda)}$ and $\mathcal{R}_{VW}=1$. We get
$$\sum_\nu (\frac{f(\lambda+\nu)}{f(\lambda)})_{|W[\nu]}\mathrm{dim}\;W[\nu]
\frac{q^{2(\mu+\rho,\lambda+\nu)}}{\delta_q^T(\lambda+\nu)}=\chi_W(q^{2(\mu+
\rho)})\frac{q^{2(\mu+\rho,\lambda)}}{\delta^T_q(\lambda)}.$$
As in \cite{EV4}, Corollary 2.16 we conclude that one can take $f(\lambda)=\delta^T_q(\lambda)$. 
\paragraph{}Theorem~\ref{T:01} now follows from (\ref{E:06}), 
Proposition~\ref{P:02} ii), Proposition~\ref{P:03} and from the following
 easily checked fusion identity~:
\begin{equation}\label{E:23}
\mathbb{J}_T^{1\cdots N}(\lambda)^{-1}\big(\mathbb{R}_T^{0,1\cdots N}
\big)\mathbb{J}_T^{1\cdots N}(\lambda+h^{(0)})=(\mathbb{R}_T^{01}
(\lambda+h^{(2\cdots N)}))\cdots(\mathbb{R}_T^{0N}(\lambda)).
\end{equation}
\section{The twisted qKZB equations} 
\paragraph{}We will first prove that the twisted qKZB equations hold for two
 finite-dimensional $U_q(\g)$-modules $V$ and $W$. As in the preceding
 section, we start with several preliminary lemmas.
\paragraph{}Consider the following function with values in $W \otimes V
 \otimes V^* \otimes W^* \otimes U_q(\g)$, with components labeled as
 $1,2,*2,*1$ and $3$ respectively~:
 $$Z_{WV}(\lambda,\mu)=\mathrm{Tr}_{|M_\mu}(\Phi^W_{\mu'+h^{(*2)}}
 \mathcal{R}^{30} \Phi^V_{\mu'}q_0^{2\lambda}B_0).$$
\begin{lem}\label{L:10}
We have
\begin{equation}\label{E:09}
Z_{WV}(\lambda,\mu)=(\mathcal{R}^{32})^{-1}\mathcal{J}_T^{12,3}(\lambda)
 \Psi_{WV}^T(\lambda+\frac{1}{2}h^{(3)},\mu).
\end{equation}
\end{lem}
\noindent
\textit{Proof.} Moving the R-matrix around the trace and using the cyclicity
 property, we get
\begin{equation*}
\begin{split}
Z_{WV}(\lambda,\mu)&=\mathcal{R}^{31}\mathrm{Tr}_{|M_\mu}(\Phi^W_{\mu'+
h^{(*2)}}\Phi^V_{\mu'}q_0^{2\lambda}B_0 \mathcal{R}^{30})\\
&=\mathcal{R}^{31}q_{12}^{2\lambda}B_3^{-1}\mathrm{Tr}_{|M_\mu}(\Phi^W_{\mu'+
h^{(*2)}} \mathcal{R}^{30} \Phi^V_{\mu'}q_0^{2\lambda}B_0)\\
&=\mathcal{R}^{31}q_{12}^{2\lambda}(B_3^{-1}\mathcal{R}^{32})B_3^{-1}Z_{WV}
(\lambda,\mu).
\end{split}
\end{equation*}
On the other hand, 
\begin{equation*}
\begin{split}
\mathcal{R}^{31}q_{12}^{2\lambda}(B_3^{-1}\mathcal{R}^{32})B_3^{-1}
\big[(\mathcal{R}^{32})^{-1}\mathcal{J}_T^{12,3}&(\lambda) \Psi_{WV}^T(\lambda
+\frac{1}{2}h^{(3)},\mu)\big]\\
&=\mathcal{R}^{31}q_{12}^{2\lambda}B_{12}\mathcal{J}_T^{12,3}(\lambda) 
\Psi_{WV}^T(\lambda+\frac{1}{2}h^{(3)},\mu).
\end{split}
\end{equation*}
From the modified ABRR equation it follows that
$$\mathcal{R}^{3,12}q_{12}^{2\lambda}B_{12}(\mathcal{J}_T^{12,3}(\lambda))=
\mathcal{J}_T^{12,3}(\lambda)q_{13}^{\Omega_\l}q_{23}^{\Omega_\l}
q_{12}^{2\lambda}.$$
Using the coproduct formula $\mathcal{R}^{3,12}=\mathcal{R}^{32}
\mathcal{R}^{31}$ and the $\l$-invariance of $Z_{VW}(\lambda,\mu)$, we see that
\begin{equation*}
\begin{split}
\mathcal{R}^{31}q_{12}^{2\lambda}(B_3^{-1}\mathcal{R}^{32})B_3^{-1}
\big[(\mathcal{R}^{32})^{-1}\mathcal{J}_T^{12,3}(\lambda) &\Psi_{WV}^T
(\lambda+\frac{1}{2}h^{(3)},\mu)\big]\\
&=(\mathcal{R}^{32})^{-1}\mathcal{J}_T^{12,3}(\lambda) 
\Psi_{WV}^T(\lambda+\frac{1}{2}h^{(3)},\mu).
\end{split}
\end{equation*}
Thus both sides of (\ref{E:09}) satisfy the equation $X=\mathcal{R}^{31}
q_{12}^{2\lambda}B_3^{-1}X$ and are of the form $X=q_{32}^{\Omega_\h}
q^{Z}_{12,3}\Psi^T_{VW}(\lambda+\frac{1}{2}h^{(3)},\mu)+\;l.o.t$. But it is
 easy to see that such an $X$ is unique and the lemma follows.\qed
\paragraph{}Now set
$$\tilde{Z}_{WV}(\lambda,\mu)=m_{32}S_3(Z_{WV}(\lambda,\mu))=
\mathrm{Tr}_{|M_\mu}(\Phi^W_{\mu'+h^{*2}}(\mathcal{R}^{20})^{-1}\Phi_{\mu}^V
 q_0^{2\lambda}B_0).$$
\begin{lem}We have
\begin{equation}\label{E:10}
\tilde{Z}_{WV}(\lambda,\mu)=S(u)^{-1}_2\mathbb{Q}_T(\lambda-\frac{1}{2}h^{(1)})_2
\mathcal{J}_T^{-1}(\lambda-\frac{1}{2}\hat{h}^{(2)})\Psi^T_{WV}(\lambda-
\frac{1}{2}\hat{h}^{(2)},\mu).
\end{equation}
\end{lem}
\noindent
\textit{Proof.} From (\ref{E:09}) it follows that
\begin{equation*}
\begin{split}
\tilde{Z}_{VW}(\lambda,\mu)&=m_{32}S_3\big[(\mathcal{R}^{32})^{-1}
\mathcal{J}_T^{12,3}(\lambda)\Psi_{WV}^T(\lambda+\frac{1}{2}h^{(3)})\big]\\
&=\sum_{ij} c_j^{(1)} \otimes S(d_j(\lambda))S(a'_i)b'_ic_j^{(2)}\Psi^T_{WV}
(\lambda-\frac{1}{2}\hat{h}^{(2)},\mu).
\end{split}
\end{equation*}
To conclude the proof of the lemma we use the following relations~: $\sum_i
 S(a'_i)b'_i=S(u^{-1})$, $S(u^{-1})x=S^2(x)S(u^{-1})$ for all $x \in U_q(\g)$
 and
$$\sum_j c^{(1)}_j\otimes S^{-1}(d_j(\lambda))c_j^{(2)}=\mathbb{Q}_T(\lambda-
\frac{1}{2}
h^{(1)})_2\mathcal{J}_T^{-1}(\lambda-\frac{1}{2}\hat{h}^{(2)}).$$
This last equation is obtained by applying $m_{32}(1 \otimes 1 \otimes
 S^{-1})$ to the cocycle identity (\ref{E:2cocy}).\qed

\paragraph{}Consider the following function with values in $V \otimes V^*
 \otimes U_q(\g)$~:
$$Y_V(\lambda,\mu)=\mathrm{Tr}_{|M_\mu}(\Phi^V_{\mu'}(\mathcal{R}^{02})^{-1}
q_0^{2\lambda}B_0).$$
\begin{lem}We have
\begin{equation}\label{E:11}
Y_V(\lambda,\mu)=q_{1}^{-2\lambda}\mathcal{J}_T^{21}(\lambda)\Psi^T_V(\lambda-
\frac{1}{2}h^{(2)},\mu).
\end{equation}
\end{lem}
\noindent
\textit{Proof.} A computation similar to the one in Lemma~\ref{L:10} shows
 that both $Y_V(\lambda,\mu)$ and $q_1^{-2\lambda}\mathcal{J}_T^{21}(\lambda)
\Psi^T_V(\lambda-\frac{1}{2}h^{(2)},\mu)$ satisfy the equation
$$X=(\mathcal{R}^{12})^{-1}q_1^{2\lambda}B_2^{-1}X$$
and are of the form $X=q_1^{-2\lambda}q_{12}^Z \Psi^T_V(\lambda-\frac{1}{2}
h^{(2)},\mu)+\;l.o.t$. It is easy to see that such an $X$ is unique.\qed
\paragraph{}We will also need the following two-representations analogue of
 $Y_V(\lambda,\mu)$~:
$$Y_{WV}(\lambda,\mu)=\mathrm{Tr}_{|M_\mu}(\Phi^W_{\mu'+h^{(*2)}}
\Phi^V_{\mu'}(\mathcal{R}^{03})^{-1}q_0^{2\lambda}B_0).$$
\begin{lem}We have
\begin{equation}\label{E:12}
Y_{WV}(\lambda,\mu)=(\mathcal{R}^{12,3})^{-1}\mathcal{J}_T^{3,12}(\lambda)
\Psi^T_{WV}(\lambda-\frac{1}{2}h^{(3)},\mu).
\end{equation}
\end{lem}
\noindent
\textit{Proof.} One checks that both sides of (\ref{E:12}) are solutions of
 the equation
$$X=(\mathcal{R}^{12,3})^{-1}q_{12}^{2\lambda}B_3^{-1}X$$ of the form
 $X=q_{12}^{-\Omega_\h}q_{13}^{-\Omega_\h}q_{31}^Zq_{32}^Z\Psi_{WV}^T(\lambda-
\frac{1}{2}h^{(3)})+\;l.o.t$.\qed
\paragraph{}Finally, we introduce a function~:
$$\tilde{Y}_{WV}(\lambda,\mu)=m_{32}(q_2^{2\rho}S_3(Y_{WV}(\lambda,\mu)))=
\mathrm{Tr}_{|M_\mu}(\Phi^W_{\mu'+h^{(*2)}}m_{32}(q_2^{2\rho}\Phi^V_{\mu'}
\mathcal{R}^{03})q_0^{2\lambda}B_0).$$
\begin{lem}\label{L:12}
We have
\begin{equation}\label{E:13}
\tilde{Y}_{WV}(\lambda,\mu)=(q^{2\rho}u^{-1}S(\mathbb{Q}_T)(\lambda-\frac{1}{2}h))_2(
\mathbb{R}^{21}\mathcal{J}_T^{-1}\Psi_{WV}^T)(\lambda+\frac{1}{2}\hat{h}^{(2)}
,\mu).
\end{equation}
\end{lem}
\noindent
\textit{Proof.} By (\ref{E:12}) we have
$$\tilde{Y}_{WV}(\lambda,\mu)=\big(\sum_{ij}(a'_i)^{(1)}d_j^{(1)}(\lambda) 
\otimes S(c_j)S(b'_i)q^{2\rho}(a'_i)^{(2)}d_j^{(2)}(\lambda)\big)\Psi_{WV}^T
(\lambda+\frac{1}{2}\hat{h}^{(2)},\mu).$$
Now we use the following relations: $q^{2\rho}x=S^2(x)q^{2\rho}$ for any $x
 \in U_q(\g)$, 
\begin{equation}\label{E:15}
\sum_i (a'_i)^{(1)}\otimes S^{-1}(b'_i)(a'_i)^{(2)}=\sum_i a'_i \otimes
 S^{-1}(b'_i)u^{-1}
\end{equation}
and $(1 \otimes S)\mathcal{R}^{-1}=\mathcal{R}$. We get
$$\tilde{Y}_{WV}(\lambda,\mu)=(q^{2\rho}u^{-1})_2\big(\sum_{ij}a_id_j^{(1)}
(\lambda)\otimes S(c_j)b_id_j^{(2)}(\lambda)\big)\Psi^T_{WV}(\lambda+
\frac{1}{2}\hat{h}^{(2)},\mu).$$
Using the identity $\mathcal{R}^{12}\mathcal{J}_T^{3,12}(\lambda)=
\mathcal{J}_T^{3,21}(\lambda)\mathcal{R}^{12}$ and
$$\sum_i S(c_i)d_i^{(1)}(\lambda) \otimes d_i^{(2)}(\lambda)=S(\mathbb{Q}_T)(\lambda-
\frac{1}{2}h^{(2)})_1\mathcal{J}_T^{-1}(\lambda-\frac{1}{2}\hat{h}^{(1)})$$
which is obtained by applying $m_{12}(S \otimes 1 \otimes 1)$ to 
(\ref{E:2cocy}), we can further simplify $\tilde{Y}_{WV}(\lambda,\mu)$ :
\begin{equation*}
\begin{split}
&\tilde{Y}_{WV}(\lambda,\mu)\\
&=(q^{2\rho}u^{-1})_2S(\mathbb{Q}_T)(\lambda-\frac{1}{2}h^{(1)})_2 
(\mathcal{J}_T^{21})^{-1}(\lambda+\frac{1}{2}\hat{h}^{(2)})\mathcal{R}^{12}
\Psi_{WV}^T(\lambda+\frac{1}{2}\hat{h}^{(2)})\\
&=(q^{2\rho}u^{-1})_2S(\mathbb{Q}_T)(\lambda-\frac{1}{2}h^{(1)})_2 \mathbb{R}^{21}_T
(\lambda+\frac{1}{2}\hat{h}^{(2)})\mathcal{J}_T^{-1}(\lambda+\frac{1}{2}
\hat{h}^{(2)})\Psi_{WV}^T(\lambda+\frac{1}{2}\hat{h}^{(2)},\mu).
\end{split}
\end{equation*}
which proves the Lemma. \qed
\paragraph{}We need one last technical result :
\begin{lem}We have
\begin{equation}\label{E:14}
m_{21}(q^{2\rho}_1\Phi^V_{\mu'}\mathcal{R}^{02})=q_{*1}^{-2(\mu'+\rho)-\sum_i 
x_i^2}(\mathcal{R}^{10})^{-1}\Phi_{\mu'}^V
\end{equation}
where $(x_i)_i$ is an orthonormal basis of $\h$.\end{lem}
\noindent
\textit{Proof.} Let $Z=u^{-1}q^{2\rho}$. This is a ribbon element of $U_q(\g)$ (see \cite{D}). 
Thus
$$\Phi^V_{\mu'}Z=\Delta(Z)\Phi^V_{\mu'}=\mathcal{R}^{21}\mathcal{R}(Z 
\otimes Z)\Phi^V_{\mu'}.$$
But by (\ref{E:15}) we have
\begin{equation*}
\begin{split}
m_{21}(q_1^{2\rho}\Phi^V_{\mu'}\mathcal{R}^{02})=&m_{21}(q_1^{2\rho}\mathcal{R}^{02}\mathcal{R}^{12})\Phi^V_{\mu'}\\
=&\mathcal{R}^{01}m_{21}(q_1^{2\rho}\mathcal{R}q_1^{-2\rho})q_1^{2\rho}\Phi^V_{\mu'}\\
=&\mathcal{R}(1 \otimes Z)\Phi^V_{\mu'}.
\end{split}
\end{equation*}
On the other hand, it is easy to see that $Z_{|M_\nu}=q^{(2\rho+\nu,\nu)}$. The Lemma now follows by a direct computation.\qed
\paragraph{}We are now in position to prove Theorem~\ref{T:02}. From
 (\ref{E:14}), we have
\begin{equation*}
\begin{split}
(q^{-2(\mu'+\rho)-\sum_i x_i^2})_{*2}\mathrm{Tr}_{|M_\mu}(\Phi^W_{\mu'+h^{*2}}
(&\mathcal{R}^{20})^{-1}\Phi^V_{\mu'}q_0^{2\lambda}B_0)\\
&=\mathrm{Tr}_{|M_{\mu}}(\Phi^W_{\mu'+h^{*2}}m_{32}(q_2^{2\rho}\Phi^V_{\mu'}
\mathcal{R}^{03})q_0^{2\lambda}B_0).
\end{split}
\end{equation*}
In other words,
$$(q^{-2(\mu'+\rho)-\sum_i x_i^2})_{*2}\tilde{Z}_{WV}(\lambda,\mu)=
\tilde{Y}_{WV}(\lambda,\mu).$$
Using (\ref{E:10}), (\ref{E:13}), the relation $uS(u^{-1})=q^{4\rho}$, 
Proposition~\ref{P:03} and the definition of $\varphi^T_{WV}(\lambda)$ we
 finally obtain
\begin{equation}\label{E:16}
(q^{-2(\mu'+\rho)-\sum_{i}x_i^2})_{*2}\varphi^T_{WV}(\lambda-h^{(2)},\mu)=
\mathbb{R}^{21}_T(\lambda)\varphi^T_{WV}(\lambda,\mu).
\end{equation}
%Using the relation $\mu'=\mu-h^{(*1)}-h^{(*2)}$ it is easy to deduce from (\ref{E:16}) the twisted qKZB equation for two representations.
From this we 
derive the qKZB equation with N representations in the following way. We
 start with the easily checked fusion identities
\begin{align}
\mathbb{J}_T^{23}(\lambda)^{-1}\mathbb{R}_T^{1,23}(\lambda)\mathbb{J}_T^{23}
(\lambda+h^{(1)})&=\mathbb{R}_T^{12}(\lambda+h^{(3)})\mathbb{R}_T^{13}
(\lambda),\label{F:01}\\
\mathbb{J}_T^{12}(\lambda+h^{(3)})^{-1}\mathbb{R}_T^{12,3}(\lambda)
\mathbb{J}_T^{12}(\lambda)&=\mathbb{R}_T^{23}(\lambda)\mathbb{R}_T^{13}
(\lambda+h^{(2)}).\label{F:02}
\end{align}
Now, from (\ref{E:16}) with $W=V_1\otimes \cdots \otimes V_j$ and $V=V_{j+1} \otimes \cdots \otimes V_N$ we get
\begin{equation}\label{F:03}
\begin{split}
(&q^{-2(\mu+\rho)+\sum x_i^2})_{*j+1\ldots, *N}(q^{2\sum x_i \otimes x_i})_{*j+1,\ldots, *N,*1,\ldots, *j}\mathbb{J}_T^{1\ldots, j}
(\lambda)\mathbb{J}_T^{j+1\ldots N}(\lambda-h^{(j+1\ldots N)})\\
&\;\times
\varphi^T_{V_1,\ldots V_N}(\lambda-h^{(j+1\ldots N)},\mu)\\
&=\mathbb{R}_T^{j+1\ldots N,1\ldots j}(\lambda) \mathbb{J}_T^{1\ldots j}
(\lambda+h^{(j+1\ldots N)})\mathbb{J}_T^{j+1\ldots N}(\lambda)
\varphi^T_{V_1,\ldots V_N}(\lambda,\mu).
\end{split}
\end{equation}
By (\ref{F:01}) this implies that
\begin{equation}\label{F:04}
\begin{split}
(&q^{-2(\mu+\rho)+\sum x_i^2})_{*j+1\ldots *N}(q^{2\sum x_i \otimes x_i})_{*j+1,\ldots *N,*1,\ldots *j}\mathbb{J}_T^{1\ldots j-1}
(\lambda+h^{(j)})\mathbb{J}_T^{j+1\ldots N}(\lambda-h^{(j+1\ldots N)})\\
&\;\times\varphi^T_{V_1,\ldots, V_N}(\lambda-h^{(j+1\ldots N)},\mu)\\
&=\mathbb{R}_T^{j+1\ldots N,1\ldots j-1}(\lambda+h^{(j)})
\mathbb{R}_T^{j+1\ldots N, j}(\lambda) \mathbb{J}_T^{1\ldots j-1}
(\lambda+h^{(j\ldots N)})\mathbb{J}_T^{j+1\ldots N}(\lambda)
\varphi^T_{V_1,\ldots, V_N}(\lambda,\mu).
\end{split}
\end{equation}
On the other hand, by (\ref{F:01}) with $W=V_1 \otimes \cdots \otimes V_{j-1}$ and $V=V_j \otimes \cdots \otimes V_N$ and using (\ref{F:02}) we also have
\begin{equation}\label{F:05}
\begin{split}
(&q^{-2(\mu+\rho)+\sum x_i^2})_{*j\ldots *N}(q^{2\sum x_i \otimes x_i})_{*j,\ldots *N,*1,\ldots *j-1}\mathbb{J}_T^{1\ldots j-1}
(\lambda)\mathbb{J}_T^{j+1\ldots N}(\lambda+h^{(1\ldots j-1)})\\
&\;\times\varphi^T_{V_1,\ldots, V_N}(\lambda-h^{(j\ldots N)},\mu)\\
&=\mathbb{R}_T^{j+1\ldots N,1\ldots j-1}(\lambda) \mathbb{R}_T^{j,1\ldots j-1}
(\lambda+h^{(j+1\ldots N)})\mathbb{J}_T^{1\ldots j-1}(\lambda+h^{(j\ldots N)})
\mathbb{J}_T^{j+1\ldots N}(\lambda)\\
&\;\;\times\varphi^T_{V_1,\ldots, V_N}(\lambda,\mu).
\end{split}
\end{equation}
Applying the operator $\Gamma_j$ to both sides of (\ref{F:05}) we get
\begin{equation}\label{F:06}
\begin{split}
(&q^{-2(\mu+\rho)+\sum x_i^2})_{*j\ldots *N}(q^{2\sum x_i \otimes x_i})_{*j,\ldots *N,*1,\ldots *j-1}\mathbb{J}_T^{1\ldots j-1}
(\lambda+h^{(j)})\mathbb{J}_T^{j+1\ldots N}(\lambda+h^{(1\ldots j)})\\
&\;\times \varphi^T_{V_1,\ldots, V_N}(\lambda-h^{(j+1\ldots N)},\mu)\\
&=\mathbb{R}_T^{j+1\ldots N,1\ldots j-1}(\lambda+h^{(j)}) \Gamma_j
\mathbb{R}_T^{j,1\ldots j-1}(\lambda+h^{(j+1\ldots N)})
\mathbb{J}_T^{1\ldots j-1}(\lambda+h^{(j\ldots N)})\mathbb{J}_T^{j+1\ldots N}
(\lambda)\\
&\;\;\times\varphi^T_{V_1,\ldots, V_N}(\lambda,\mu).
\end{split}
\end{equation}
Comparing (\ref{F:04}) with (\ref{F:06}) and using the $\l$-invariance of $\mathbb{R}$ we obtain
\begin{equation}\label{F:07}
\begin{split}
(&q^{-2(\mu+\rho)+\sum x_i^2})_{*j}(q^{2\sum x_i \otimes x_i}
)_{*j,*1\ldots *j-1}\mathbb{J}_T^{1\ldots j-1}(\lambda+h^{(j\ldots N)})
\mathbb{R}_T^{j+1\ldots N,j}(\lambda)\mathbb{J}_T^{j+1\ldots N}(\lambda)\\
&\;\times\varphi^T_{V_1,\ldots, V_N}(\lambda,\mu)\\
&=\mathbb{J}_T^{j+1\ldots N}(\lambda+h^{(j)})\Gamma_j 
\mathbb{R}_T^{j,1\ldots j-1}(\lambda+h^{(j+1\ldots N)})
\mathbb{J}_T^{1,\ldots j-1}(\lambda+h^{(j\ldots N)})\varphi_{V_1,\ldots, V_N}^T
(\lambda,\mu),
\end{split}
\end{equation}
which can be rewritten as
\begin{equation}\label{F:08}
\begin{split}
(&q^{-2(\mu+\rho)+\sum x_i^2})_{*j}(q^{2\sum x_i \otimes x_i}
)_{*j,*1\ldots *j-1}\mathbb{J}_T^{j+1\ldots N}(\lambda+h^{(j)})^{-1}
\mathbb{R}_T^{j+1\ldots N,j}(\lambda)\mathbb{J}_T^{j+1\ldots N}(\lambda)\\
&\;\times\varphi^T_{V_1,\ldots, V_N}(\lambda,\mu)\\
&=\Gamma_j\mathbb{J}_T^{1\ldots j-1}(\lambda+h^{(j+1\ldots N)})^{-1} 
\mathbb{R}_T^{j,1\ldots j-1}(\lambda+h^{(j+1\ldots N)})
\mathbb{J}_T^{1,\ldots j-1}(\lambda+h^{(j\ldots N)})\\
&\;\;\times\varphi_{V_1,\ldots, V_N}^T
(\lambda,\mu).
\end{split}
\end{equation}
Finally, taking into account identities (\ref{F:01}) and (\ref{F:02}), we
 obtain
\begin{equation}\label{F:09}
\begin{split}
&(q^{-2(\mu+\rho)+\sum x_i^2})_{*j}(q^{2\sum x_i \otimes x_i}
)_{*j,*1\ldots *j-1}\mathbb{R}_T^{Nj}(\lambda)\cdots \mathbb{R}_T^{j+1,j}
(\lambda+h^{(j+2\ldots N)})\varphi^T_{V_1,\ldots ,V_N}(\lambda,\mu)\\
&=\Gamma_j\mathbb{R}_T^{j1}(\lambda+h^{(2\ldots j-1)}+h^{(j+1\ldots N)})
\cdots \mathbb{R}_T^{jj-1}(\lambda+h^{(j+1\ldots N)})
\varphi^T_{V_1,\ldots, V_N}(\lambda,\mu).
\end{split}
\end{equation}
The proof of Theorem~\ref{T:02} is now obtained by replacing $\mu$ by $-\mu-\rho$ and by rewriting (\ref{F:09}) in terms of $F_{V_1,\ldots, V_N}^T(\lambda,\mu)$.
\section{The twisted dual Macdonald-Ruijsenaars equation}
\paragraph{}In this section we let $(\Gamma_1,\Gamma_2,T)$ be a complete generalized Belavin-Drinfeld triple. Let $W$ be a finite-dimensional $U_q(\g)$-module such that $W \simeq W^B$ and let us consider $W$ as a $\langle B \rangle \ltimes U_q(\g)$-module as in Section 2.
\paragraph{}For generic values of $\mu$, the tensor product $M_\mu \otimes {W}$ decomposes as a direct sum of Verma modules, and 
\begin{align*}
\eta_\nu:\; {W}[\nu] \otimes M_{\mu+\nu} &\to M_{\mu} \otimes {W}\\
w \otimes y &\mapsto \Phi^w_{\mu+\nu}(y)
\end{align*}
is an isomorphism onto the isotypic component corresponding to $M_{\mu+\nu}$. The following lemma is straightforward :
\begin{lem}\label{L:d1} We have $(B \otimes B) \circ \eta_{\nu}=\eta_{B(\nu)}\circ (B \otimes B)$.\end{lem}
Now let $V$ be any finite-dimensional $U_q(\g)$-module and consider the composition
$$P_{V \otimes V^*,{W}} \mathcal{R}^{V{W}}\Phi^V_{B(\mu)}(B \otimes B)\eta_\nu:\; {W}[\nu]\otimes M_{\mu+\nu}\to M_{B(\mu)+h^{(V)}} \otimes {W}\otimes V \otimes V^*.$$
By \cite{EV4}, Proposition 3.1, we have 
$$P_{V \otimes V^*,{W}} \mathcal{R}^{V{W}}\Phi^V_{B(\mu)}\eta_{\nu}=R^{{W}V}(B(\mu+\nu))^{t_2}\Phi^V_{B(\mu+\nu)}.$$
It follows from Lemma~\ref{L:d1} that
\begin{equation}\label{E:20}
P_{V \otimes V^*,{W}} \mathcal{R}^{V{W}}\Phi^V_{B(\mu)}(B \otimes B)\eta_\nu=\eta_{h^{({W})}}R^{{W}V}(B(\mu+\nu))^{t_2}\Phi^V_{B(\mu+\nu)} (B \otimes B)
\end{equation}
where $R(\lambda)=R_{Id}(\lambda)$ is the quantum dynamical R-matrix corresponding to the trivial triple $(\Gamma,\Gamma,Id)$ and where $t_2$ means transposition in the second component (so that $R^{{W}V}(B(\mu+\nu))^{t_2}$ acts on ${W}\otimes V^*$). Now let us multiply both sides of (\ref{E:20}) by $q^{2\lambda}_{M_\mu \otimes {W}}$ and sum over all values of $\nu$. This yields
\begin{equation}\label{E:21}
P_{V \otimes V^*,{W}}\mathcal{R}^{V{W}}\Phi^V_{B(\mu)}(B \otimes B)q^{2\lambda}_{M_\mu \otimes {W}}=\eta R^{{W}\otimes V}(B(\mu+h^{({W})})) (B \otimes B) q^{2\lambda} \eta^{-1}
\end{equation}
where $\eta=\bigoplus_\nu \eta_\nu:\; \bigoplus_\nu {W}[\nu]\otimes M_{\mu+\nu} \stackrel{\sim}{\to} M_\mu \otimes {W}$. Let us take the trace in the Verma modules and in ${W}$. Using the fact that $\mathcal{R} \in q^{\Omega_\h} U_q(\n_+) \otimes U_q(\n_-)$ and that $\nu-B(\nu)$ is never a linear combination of negative roots, we obtain
$$\mathrm{Tr}_{|{W}}(q^{2\lambda+h^{(V)}}B) \varphi^T_V(\lambda,\mu)=\sum_\nu \mathrm{Tr}_{|{W}[\nu]}(R^{{W}V}(B(\mu+\nu))^{t_2}B) \varphi^T_{V}(\lambda,\mu+\nu).$$
It is clear that
$$\mathrm{Tr}_{|{W}}(q^{2\lambda+h^{(V)}}B)=\sum_{\nu, B(\nu)=\nu} \mathrm{Tr}_{|W[\nu]}(q^{2\lambda+h^{(V)}}B)=\mathrm{Tr}_{|W^{\h_0}}(q^{2\lambda}B).$$
Hence from (\ref{E:21}) we get
$$\mathrm{Tr}_{|(W^*)^{\h_0}}(q^{-2\lambda}B^*)\varphi^T_V(\lambda,\mu)=\sum_{\nu} \mathrm{Tr}_{|W^*[-\nu]}(B^*_{W^*}R^{WV}(B(\mu+\nu))^{t_1t_2})\varphi^T_V(\lambda,\mu+\nu),$$
which can be rewritten in terms of $F^T_V(\lambda,\mu)$ as
\begin{equation}\label{E:22}
\begin{split}
&\mathrm{Tr}_{|(W^*)^{\h_0}}(q^{-2\lambda}B)F^T_V(\lambda,\mu)\\
&=\sum_{\nu \in \l^*} \mathrm{Tr}_{|W^*[\nu]}\big(\mathbb{Q}^{-1}_{|V^*}(B(\mu))B^*_{W^*}\mathbb{R}^{WV}(B(\mu+\nu))^{t_1t_2}\mathbb{Q}_{|V^*}(B(\mu+\nu))\big)F^T_V(\lambda,\mu+\nu).
\end{split}
\end{equation}
Finally, using the formula
\begin{equation}\label{E:312}
\mathbb{R}_{WV}(\lambda)^{t_1t_2}=(\mathbb{Q}(\lambda) \otimes \mathbb{Q}(\lambda-h^{(1)}))\mathbb{R}_{W^*V^*}(\lambda-h^{(1)}-h^{(2)})(\mathbb{Q}^{-1}(\lambda-h^{(2)})\otimes \mathbb{Q}^{-1}(\lambda))
\end{equation}
 (see \cite{EV4}, (3.12)) and using the fact that $\mathbb{Q}$ is of weight zero, we simplify (\ref{E:22}) to
\begin{equation}\label{E:40}
\begin{split}
\mathrm{Tr}&_{|(W^*)^{\h_0}}(q^{-2\lambda}) F_V^T(\lambda,\mu)\\
&=\sum_\nu \mathrm{Tr}_{|W^*[\nu]}\big(B^*_{W^*}\mathbb{Q}_{W^*}(B(\mu+\nu))\mathbb{R}_{W^*V^*}(B(\mu+\nu)-\nu-h^{(2)})\times\\
&\qquad \qquad \qquad \qquad \qquad \qquad \qquad \qquad \times \mathbb{Q}^{-1}_{W^*}(B(\mu+\nu)-h^{(2)})\big)F^T_V(\lambda,\mu+\nu)\\
&=\sum_\nu \mathrm{Tr}_{|W^*[\nu]}\big(\mathbb{Q}_{W^*}(B(\mu+\nu))\mathbb{R}_{W^*V^*}(\mu)B^*_{W^*}\mathbb{Q}^{-1}_{W^*}(B(\mu+\nu))\big)F^T_V(\lambda,\mu+\nu)\\
&=\sum_\nu \mathrm{Tr}_{|W^*[\nu]}\big(\mathbb{R}_{W^*V^*}(\mu)B^*_{W^*}\big)F^T_V(\lambda,\mu+\nu).
\end{split}
\end{equation}
The twisted dual Macdonald-Ruijsenaars equations for an arbitrary number of modules $V_1,\ldots V_N$ can now be deduced from (\ref{E:40}) and from the fusion identity (\ref{E:23}). Theorem~\ref{T:03} is proved.
\noindent
\paragraph{}\textit{Proof of Proposition~\ref{P:001}.} Let $\mathbf{W}$ and $\overline{\mathbf{W}}$ be the Weyl groups of $\Gamma$ and $\overline{\Gamma}$ respectively. By the Bernstein-Gelfand-Gelfand resolution, we have
$$\mathrm{Tr}_{|V_\nu}(q^{2\lambda}B)=\sum_{w \in \mathbf{W}} (-1)^{l(w)}\mathrm{Tr}_{|M_{w(\nu+\rho)-\rho}}(q^{2\lambda}B).$$
Denote by $s_\alpha$ the simple reflection corresponding to the simple root $\alpha \in \Gamma$. The group generated by $B$ acts on $\mathbf{W}$ by $B(s_\alpha)=s_{T\alpha}$. It follows from the facts that $\mathbf{W}$ acts simply transitively on the sets of simple roots and that $\nu$ is dominant that $B(w(\nu+\rho)-\rho)=w(\nu+\rho)-\rho$ if and only if $B(w)=w$. Moreover, $\mathbf{W}^B$ is naturally isomorphic to $\overline{\mathbf{W}}$. Hence,
\begin{equation}
\begin{split}
\sum_{w \in W} (-1)^{l(w)}\mathrm{Tr}_{|M_{w(\nu+\rho)-\rho}}(q^{2\lambda}B)&=\sum_{w \in W^B} (-1)^{l(w)}\mathrm{Tr}_{|M_{w(\nu+\rho)-\rho}}(q^{2\lambda}B)\\
&=\sum_{w\in \overline{W}}(-1)^{l(w)}\frac{q^{2(\lambda,w(\nu+\rho)-\rho)}}{\prod_{\overline{\alpha}\in \overline{\Delta}^+} (1-\theta_{\alpha}q^{-2(\overline{\alpha},\lambda)})}
\end{split}
\end{equation}
 Let $\omega_\alpha$ be the fundamental weight corresponding to $\alpha \in \Gamma$. It is easy to check that $\{\overline{\omega}_{\overline{\alpha}}=\frac{(\overline{\alpha},\overline{\alpha})}{2N_\alpha}\sum_{\alpha \in \overline{\alpha}} \omega_{{\alpha}}, \overline{\alpha} \in \overline{\Gamma}\}$ is the set of fundamental weights of $\overline{\Gamma}$. Thus $(2\lambda, w(\nu+\rho)-\rho)=(2\overline{\lambda},w(\nu+\overline{\rho})-\overline{\rho})$ where $\overline{\rho}=\sum_{\overline{\alpha}} \overline{\omega}_{\overline{\alpha}}$. Hence, by the Weyl character formula for $\overline{\g}$ we have
$$\mathrm{Tr}_{|V_\nu}(q^{2\lambda}B)=\chi_{\overline{V}_\nu}(q^{2\overline{\lambda}})\frac{\overline{\delta}_q(\overline{\lambda})}{\delta_q^T(\lambda)}
$$
where
 $$\overline{\delta}_q(\overline{\lambda})=q^{2(\rho,\overline{\lambda})}\prod_{\overline{\alpha} \in \overline{\Delta}^+}(1-q^{-2(\overline{\alpha},\overline{\lambda})})$$
is the (usual) Weyl denominator for $\overline{\Gamma}$. Setting $\nu=0$ we see that in fact $\overline{\delta}_q(\overline{\lambda})= \delta_q^T(\lambda)$ . The Proposition follows.\qed
\section{The twisted dual qKZB equations}
\paragraph{}In this section we prove Theorem~\ref{T:04}. As in the preceding section, let $T$ be an automorphism of $\Gamma$ and let $V_1,\ldots V_N$ be finite-dimensional $U_q(\g)$-modules. We will extensively use the following two identities which are proved in \cite{EV3} :
\begin{equation}\label{E:50}
\Phi_{\mu+h^{(V^*)}}^{W}\Phi_{\mu}^{V}=\mathcal{R}^{-1}
R_{21}^*(\mu)\Phi_{\mu+h^{(W^*)}}^{V}\Phi_{\mu}^{W}=
\mathcal{R}_{21}
R^*(\mu)^{-1}\Phi_{\mu+h^{(W^*)}}^{V}\Phi_{\mu}^{W}
\end{equation}
for any two modules $V,W$. 
\paragraph{}Consider 
$$\Psi^T_{V_1,\ldots, V_N}(\lambda,\mu)=\mathrm{Tr}_{|M_\mu}\big(\Phi^{V_1}_{B(\mu)+h^{(*2\cdots *N)}}\cdots \Phi^{V_N}_{B(\mu)}q^{2\lambda}B\big)$$
and move the $j$th intertwiner to the right using (\ref{E:50}). We get
\begin{equation}\label{E:51}
\begin{split}
\Psi&^T_{V_1,\ldots, V_N}(\lambda,\mu)\\
&=\mathcal{R}_{j+1,j}\cdots \mathcal{R}_{N,j}q_j^{2\lambda}R^*_{j,j+1}(B(\mu)+h^{(*j+2\cdots *N)})^{-1}\cdots R^*_{j,N}(B(\mu))^{-1}\times\\
& \times\mathrm{Tr}_{|M_\mu}\big(\Phi^{V_1}_{B(\mu)+h^{(*2\cdots *N)}}\cdots \Phi^{V_N}_{B(\mu)+h^{(*j)}}q^{2\lambda}\Phi^{V_j}_{B(\mu)}B\big).
\end{split}
\end{equation}
Now, we have
\begin{equation}\label{E:52}
\begin{split}
\mathrm{Tr}_{|M_\mu}\big(&\Phi^{V_1}_{B(\mu)+h^{(*2\cdots *N)}}\cdots \Phi^{V_N}_{B(\mu)+h^{(*j)}}q^{2\lambda}\Phi^{V_j}_{B(\mu)}B\big)\\
&=B_{V'_j}B^*_{V^{*'}_j}\mathrm{Tr}_{M_{\mu+h^{(*j)}}}\big(\Phi^{V'_j}_{B(\mu+h^{(*j)})+\sum_{i=1,i\neq j}^N h^{(*)i}} \cdots \Phi^{V_N}_{B(\mu+h^{(*j)})}q^{2\lambda}B\big)\\
&=B_{V'_j}B^*_{V^{*'}_j}\Gamma^*_{*j}\mathrm{Tr}_{M_{\mu}}\big(\Phi^{V'_j}_{B(\mu)+\sum_{i=1,i\neq j}^N h^{(*i)}} \cdots \Phi^{V_N}_{B(\mu)}q^{2\lambda}B\big)
\end{split}
\end{equation}
where we note $V_j'=V^{B^{-1}}_j$ for simplicity.
Finally, we move $\Phi^{V'_j}$ to the right back to its original position, thereby completing a cycle. By (\ref{E:50}) we obtain
\begin{equation}\label{E:53}
\begin{split}
\mathrm{Tr}_{|M_{\mu}}&\big(\Phi^{V'_j}_{B(\mu)+\sum_{i=1,i\neq j}^N h^{(*i)}} \cdots \Phi^{V_N}_{B(\mu)}q^{2\lambda}B\big)\\
=&\mathcal{R}^{-1}_{j,1}\cdots \mathcal{R}^{-1}_{j,j-1}R^{*}_{1,j}\big(B(\mu)+\sum_{i=2,i\neq j}^N h^{(*i)}\big)\cdots R^*_{j-1,j}\big(B(\mu)+\sum_{i=j+1}^N h^{(*i)}\big) \\
&\times\Psi^T_{V_1,\ldots V'_j, \ldots V_N}(\lambda,\mu).
\end{split}
\end{equation}
Combining (\ref{E:51}), (\ref{E:52}) and (\ref{E:53}) yields the following relation
\begin{equation}\label{E:54}
\begin{split}
\Psi&^T_{V_1,\ldots, V_N}(\lambda,\mu)\\
=&\big[\mathcal{R}_{j+1,j}\cdots \mathcal{R}_{N,j}q_j^{2\lambda}(B_j\mathcal{R}_{j,1}^{-1})\cdots (B_j\mathcal{R}_{j,j-1}^{-1})\big]\times\\
&\times\big[ R^*_{j,j+1}\big(B(\mu)+\sum_{i=j+2}^N h^{(*i)}\big)^{-1}\cdots R^*_{j,N}\big(B(\mu)\big)^{-1}\Gamma^*_{B^{-1}(j)}\times\\
&\times B^*_jR^*_{1,j}\big(B(\mu)+\sum_{i=2,i\neq j}^N h^{(*i)}\big)\cdots B^*_jR^*_{j-1,j}\big(B(\mu)+\sum_{i=j+1}^N h^{(*i)}\big)\big]\times\\
&\times B_{V'_j}B^*_{V^{*'}_j} \Psi^T_{V_1,\ldots, V'_j,\ldots V_N}(\lambda,\mu).
\end{split}
\end{equation}
Let us replace $\mu$ by $-\mu-\rho$ and let us rewrite this equation in terms of $F^T(\lambda,\mu)$. We get
\begin{equation}\label{E:55}
\begin{split}
F&^T_{V_1,\ldots, V_N}(\lambda,\mu)\\
=&\big[\mathbb{J}^{1\cdots N}_T(\lambda)^{-1}\mathcal{R}_{j+1,j}\cdots \mathcal{R}_{N,j}q_j^{2\lambda}(B_j\mathcal{R}_{j,1}^{-1})\cdots (B_j\mathcal{R}_{j,j-1}^{-1})(B_j\mathbb{J}_T^{1\cdots N}(\lambda))\big]\times\\
&\times\bigg[ \bigg\{\mathbb{Q}^{-1}_{*N}(B(\mu)) \otimes \cdots \otimes \mathbb{Q}^{-1}_{*1}\big(B(\mu)-\sum_{i=2}^N h^{(*i)}\big)\bigg\}\mathbb{R}^*_{j,j+1}\big(B(\mu)-\sum_{i=j+2}^N h^{(*i)}\big)^{-1}\times\cdots\\
&\times \mathbb{R}^*_{j,N}(B(\mu))^{-1}\Gamma^{*-1}_{B^{-1}(j)}B^*_j\mathbb{R}^*_{1,j}\big(B(\mu)-\sum_{i=2,i\neq j}^N h^{(*i)}\big)\cdots B^*_j\mathbb{R}^*_{j-1,j}\big(B(\mu)-\sum_{i=j+1}^N h^{(*i)}\big)\times\\
&\times B_j\bigg\{\mathbb{Q}_{*N}(B(\mu))\otimes \cdots \otimes \mathbb{Q}_{*1}\big(B(\mu)-\sum_{i=2,i\neq j}^N h^{(*i)}-B^{-1}(h^{(*j)})\big)\bigg\}\bigg]\times\\
&\times B_{V'_j}B^*_{V^{*'}_j} F^T_{V_1,\ldots,V'_j,\ldots, V_N}(\lambda,\mu).
\end{split}
\end{equation}
Inverting, we obtain
\begin{equation}\label{E:56}
\begin{split}
B&_{V'_j}B^*_{V^{*'}_j} F^T_{V_1,\ldots, V'_j,\ldots, V_N}(\lambda,\mu)\\
=&\big[(B_j\mathbb{J}_T^{1\cdots N}(\lambda))^{-1}(B_j\mathcal{R}_{j,1\ldots j-1}) q_j^{-2\lambda} \mathcal{R}^{-1}_{j+1 \ldots N,j} \mathbb{J}_T^{1\cdots N}(\lambda)\big]\times\\
&\times B_j\bigg\{\mathbb{Q}^{-1}_{*N}(B(\mu))\otimes \cdots \otimes \mathbb{Q}^{-1}_{*1}\big(B(\mu)-\sum_{i=2,i\neq j}^N h^{(*i)}-B^{-1}(h^{(*j)})\big)\bigg\} \times\\
&\times B^*_j\mathbb{R}^{*-1}_{j-1,j}\big(B(\mu)-\sum_{i=j+1}^N h^{(*i)}\big)\cdots
\times B^*_j\mathbb{R}^{*-1}_{1,j}\big(B(\mu)-\sum_{i=2,i\neq j}^N h^{(*i)}\big)\times\\
&\times \Gamma^{*}_{B^{-1}(j)}\mathbb{R}^*_{j,N}(B(\mu))\cdots\mathbb{R}^*_{j,j+1}\big(B(\mu)-\sum_{i=j+2}^N h^{(*i)}\big)\times\\
&\times\bigg\{\mathbb{Q}_{*N}(B(\mu)) \otimes \cdots \otimes \mathbb{Q}_{*1}\big(B(\mu)-\sum_{i=2}^N h^{(*i)}\big)\bigg\}\times F^T_{V_1,\ldots, V_N}(\lambda,\mu).
\end{split}
\end{equation}
Using (\ref{E:312}) and using the fact that $\mu=B(\mu)-\sum h^{(*i)}$ it is easy to check that
\begin{equation*}
\begin{split}
K^{\vee,T}_j&=\\
&B_j\bigg\{\mathbb{Q}^{-1}_{*N}(B(\mu))\otimes \cdots \otimes \mathbb{Q}^{-1}_{*1}\big(B(\mu)-\sum_{i=2,i\neq j}^N h^{(*i)}-B^{-1}(h^{(*j)})\big)\bigg\}\times\\
&\times B^*_j\mathbb{R}^{*-1}_{j-1,j}\big(B(\mu)-\sum_{i=j+1}^N h^{(*i)}\big)\times \cdots\times B^*_j\mathbb{R}^{*-1}_{1,j}\big(B(\mu)-\sum_{i=2,i\neq j}^N h^{(*i)}\big)\times\\
&\Gamma^{*}_{B^{-1}(j)}\mathbb{R}^*_{j,N}(B(\mu))\cdots \mathbb{R}^*_{j,j+1}\big(B(\mu)-\sum_{i=j+2}^N h^{(*i)}\big)\times\\
&\times \bigg\{\mathbb{Q}_{*N}(B(\mu)) \otimes \cdots \otimes \mathbb{Q}_{*1}\big(B(\mu)-\sum_{i=2}^N h^{(*i)}\big)\bigg\}.
\end{split}
\end{equation*}
Finally, we have
\begin{equation}\label{E:last}
D^{\vee,T}_j=
(B_j\mathbb{J}_T^{1\cdots N}(\lambda))^{-1}(B_j\mathcal{R}_{j,1\ldots j-1}) q_j^{-2\lambda} \mathcal{R}^{-1}_{j+1 \ldots N,j} \mathbb{J}_T^{1\cdots N}(\lambda)
\end{equation}
when applied to $(V_1\otimes \cdots \otimes V_N)^\l$. The proof of this last equality is similar to the proof of \cite{EV4} (4.10): it is enough to check (\ref{E:last}) for $N=3$, for which it follows from the modified ABRR equation (\ref{E:MABRR}). This concludes the proof of Theorem~\ref{T:04}.
\section{The classical limits}
\paragraph{}Let us now examine the classical limits of Theorems 2.1-4, that is, the behavior of (a suitable renormalization of) the functions $F^T_{V_1,\ldots, V_N}(\lambda,\mu)$ when $\hbar \to 0$. In that limit, the quantum group $U_q(\g)$ becomes the usual enveloping algebra $U(\g)$. We will denote by $\Phi, \mathbb{Q}_T^c, \mathbb{R}_T^c, \mathbb{J}_T^c,\ldots$ the classical limits of the operators constructed in Section 2.1, obtained when we replace $U_q(\g)$ by $U(\g)$.
\paragraph{}Let $V_1,\ldots V_N$ be finite-dimensional $U_q(\g)$-modules and let $V_1^c,\ldots V_N^c$ be the corresponding $U(\g)$ modules. Let us fix a generalized Belavin-Drinfeld triple $(\Gamma_1,\Gamma_2,T)$ and set
$$\Psi^{T,c}_{V_1,\ldots, V_N}=\mathrm{Tr}_{|M^c_{\mu}}\big(\Phi^{V_1^c}_{\mu'+\sum_{i=2}^N h^{(*i)}}\cdots \Phi^{V^c_N}_{\mu'}Be^{-\lambda}\big).$$
Also set $\delta^T(\lambda)=(\mathrm{Tr}_{|M_{-\rho}^c}(B e^{-\lambda}))^{-1}$. We define the classical limit of the function $F^T_{V_1,\ldots, V_N}(\lambda,\mu)$ as
$$F^{T,c}_{V_1,\ldots, V_N}(\lambda,\mu):=\underset{\hbar \to 0}{\mathrm{lim}}\; F^T_{V_1,\ldots, V_N}(\frac{\lambda}{\hbar},\mu).$$
The following result is clear from the definitions.
\begin{lem} We have
\begin{equation*}
\begin{split}
&F^{T,c}_{V_1,\ldots, V_N}(\lambda,\mu)\\
&\quad=\delta^T(\lambda)[\mathbb{Q}^{c}_{*N}(\mu+h^{(*1\cdots *N)})^{-1}\otimes \cdots \otimes \mathbb{Q}^{c}_{*1}(\mu+h^{(*1)})^{-1}] \Psi^{T,c}_{V_1,\ldots, V_N}(\lambda,-\mu-\rho).
\end{split}
\end{equation*}
\end{lem}
\paragraph{}The classical analogue of Proposition 3.1 is as follows.
\begin{prop} Let $V$ be any finite-dimensional $U(\g)$-module and let $X \in U(\g)$.
\begin{enumerate}
\item[i)] There exists a unique differential operator $d^T_X$ acting on functions $\l^* \to V^\l$ such that
$$\mathrm{Tr}_{|M_\mu}(\Phi^V_{\mu'} X B e^{-\lambda})=d_X^T \mathrm{Tr}_{|M_\mu}(\Phi^V_{\mu'} B e^{-\lambda}).$$
\item[ii)] If $X,Y$ belong to the center of $U(\g)$ then $d_X^Td_Y^T=d_Y^Td_X^T$.
\end{enumerate}
\end{prop}
Unfortunately, there is no convenient classical analogue of the Drinfeld-Reshetikhin construction of central elements in $U_q(\g)$, and therefore no convenient explicit computation of the operator $d_X^T$ in general. However, this can be done when $X=m_{12}(\Omega)$ is the quadratic Casimir, which yields the following classical analogue of Theorem 2.1 (which is proved directly in \cite{ES1}, Theorem 3.2). 
\begin{theo}[\cite{ES1}]\label{T:71} The function
$F^{T,c}_{V_1,\ldots, V_N}(\lambda,\mu)$ satisfies the following second order differential equation :
\begin{equation}\label{E:prop12}
 \bigg(\sum_{i \in I_1} \frac{\partial^2}{\partial x_i^2}- 
\sum_{l,n=1}^r S_T(\lambda)_{|V_l \otimes V_n}\bigg)F^{T,c}_{V_1,\ldots V_N}(\lambda,\mu)
=(\mu,\mu)F^{T,c}_{V_1,\ldots, V_N}(\lambda,\mu)
\end{equation}
where $(x_i)_{i\in I_1}$ (resp. $(x_i)_{i \in I_2}$) is an orthonormal basis of $\l$ (resp. of $\h_0$) and where
\begin{equation*}
\begin{split}
S_T(\lambda)=\sum_{\alpha}\sum_{k=0}^{\infty}
\sum_{v=1}^{\infty} e^{(s+v)(\alpha,\lambda)} (B^sf_\alpha
\otimes B^{-v}e_\alpha+B^{-v}e_\alpha\otimes 
B^sf_\alpha)\\ -\sum_{i
  \in I_2} \frac{1-C_T}{2}x_i \otimes \frac{1-C_T}{2}x_i.
\end{split}
\end{equation*}
\end{theo}
\noindent
Theorem~\ref{T:71} can also be deduced from Theorem~\ref{T:01} by expanding powers of $\hbar$.
\paragraph{}The classical limit of Theorem 2.2 are the twisted (trigonometric) KZB equations.
\begin{theo}[\cite{ES1}] The function $F^{T,c}_{V_1,\ldots, V_N}(\lambda,\mu)$ satisfies the following
system of differential equations, for $j=1,\ldots N$:
\begin{equation}\label{E:prop1}
\begin{split}
\bigg(\sum_{i \in I_1} x_{i|V_j} \frac{\partial}{\partial
x_i}+&\sum_{l>j}{r}_T(\lambda)_{|V_j \otimes V_l} -
\sum_{l<j}{r}_T(\lambda)_{|V_l \otimes V_j} \bigg) F^{T,c}_{V_1,\ldots,
V_N}(\lambda,\mu)\\
&\qquad=\left((\mu+\frac{1}{2}\C_\h)_{|V^*_j}+\sum_{l=1}^{j-1}(\Omega_{\h})_{|V^*_{i} \otimes V^*_j}\right)F^{T,c}_{V_1,\ldots, V_N}(\lambda,\mu).
\end{split}
\end{equation}
\end{theo}
This theorem is proved in \cite{ES1} but can also be deduced from Theorem~\ref{T:02} by expanding in powers of $\hbar$.
\paragraph{}Finally, when $T$ is an automorphism of the Dynkin diagram $\Gamma$ we consider classical limits of the dual Macdonald-Ruijsenaars and dual qKZB equations. Let $W$ be a $B$-invariant finite-dimensional $\g$-module and let $\mathcal{D}^{\vee,T,c}_W$ denote the difference operator given by formula (\ref{E:70}) when $U_q(\g)$ is replaced by $U(\g)$.
\begin{theo} We have
$$\mathcal{D}^{\vee,T,c}_W F^{T,c}_{V_1,\ldots, V_N}=\mathrm{Tr}_{|W^{\h_0}}(e^{-\lambda}B) F^{T,c}_{V_1,\ldots, V_N}.$$
\end{theo}
Similarly, let $K^{\vee,T,c}_j$ be the classical limit of $K^{\vee,T}_j$, i.e the difference operator given by formula (\ref{E:71}) when $q=1$ (and hence $\mathbb{R}(\mu)$ is just the classical exchange matrix evaluated at $-\mu-\rho$, see \cite{EV3}).
\begin{theo} For $j=1,\ldots N$ we have
$$B_{V_j}B^*_{V_j^*}F^{T,c}_{V_1,\ldots,V_N}=(e^{-\lambda})_{|V_j}K^{\vee,T,c}_jF^{T,c}_{V_1,\ldots, V_j^B,\ldots, V_N}.$$
\end{theo}
\section{Extension to Kac-Moody algebras}
\paragraph{}In this section we briefly explain how to adapt the construction of \cite{ESS} to Kac-Moody algebras and how to generalize Theorems~\ref{T:01}-\ref{T:04} to this setting.
\paragraph{} Let $A=(a_{ij})$ be a symmetrizable
generalized Cartan matrix of size $n$ and rank $l$. Let
$(\h,\Gamma,\check{\Gamma})$ be a realization of $A$, i.e $\h$ is a complex
vector space of dimension $2n-l$, $\Gamma=\{\alpha_1,\ldots
\alpha_n\}\subset \h^*$ and $\check{\Gamma}=\{h_1,\ldots h_n\}\subset \h$
are linearly independent sets and $\langle \alpha_j,h_i\rangle =a_{ij}$.
Let $\g=\n_- \oplus \h \oplus \n_+$ be the Kac-Moody algebra associated
to $A$, i.e $\g$ is generated by elements $e_i$, $f_i$, $i=1,\ldots n$ and
$\h$ with relations
$$[e_i,f_j]=\delta_{ij}h_i,\qquad [\h,\h]=0,\qquad [h,e_i]=\langle
\alpha_i,h\rangle e_i, \qquad [h,f_i]=-\langle \alpha_i,h\rangle f_i,$$
together with the Serre relations (see \cite{K}). Let $(\,,\,)$ be a nondegenerate invariant bilinear form on $\g$. Let $\Omega_\h$ be the inverse element to the restriction $(\,,\,)$ to $\h$. For every root $\alpha \in \h^*$ we set $\alpha^\vee=(1\otimes \alpha)\Omega_\h$.
\paragraph{}Let $U_q(\g)$ be the quantum Kac-Moody algebra. It is defined by the same relations as in Section 1, where now $(a_{ij})$ is the generalized Cartan matrix $A$.
\paragraph{Construction of the twist.} Let $(\Gamma_1,\Gamma_2,T)$ be a generalized Belavin-Drinfeld triple. As before we set $\l=\big(\sum_{\alpha \in \Gamma_1}\C(\alpha -T\alpha)\big)^\perp$ and $\h_0=\l^\perp \subset \h$. We will say that $(\Gamma_1,\Gamma_2,T)$ is nondegenerate if the restriction of $(\,,\,)$ to $\l$ is, and we make this assumption from now on. Let $\h_1\subset \h$ (resp. $\h_2\subset \h$) be the subspace spanned by simple roots $\alpha \in \Gamma_1$ (resp. $\alpha \in \Gamma_2$).
\paragraph{}In \cite{ESS} we obtained an explicit construction of a twist $\mathcal{J}_T(\lambda)$ for simple complex Lie algebras. An important observation there was that $\h=\h_1+\l$, which makes it is possible to extend $B$ to an orthogonal automorphism of $\h$, and to define maps $B^{\pm 1}:\; U_q(\b_{\mp}) \to U_q(\b_{\mp})$. However, in general we only have $\h_1 +\l \subset \h$ but $\h_1+\l \neq \h$, and thus it is necessary to modify the construction in \cite{ESS}, which is done below.\\
\hbox to1em{\hfill} The following lemma is obvious.
\begin{lem}There exist unique algebra morphism $B: U_q(\n_-\oplus \h_1) \to U_q(\n_-\oplus \h_2)$ and $B: U_q(\n_+\oplus \h_2) \to U_q(\n_-\oplus \h_1)$ such that $B(F_\alpha)=F_{T\alpha}$, $B(h_\alpha)=h_{T\alpha}$ if $\alpha \in \Gamma_1$, $B(F_\alpha)=0$ if $\alpha \in \Gamma\setminus \Gamma_1$, and $B^{-1}(E_\alpha)=E_{T^{-1}\alpha}$, $B^{-1}(h_\alpha)=h_{T^{-1}\alpha}$ if $\alpha \in \Gamma_2$, $B^{-1}(E_\alpha)=0$ if $\alpha \in \Gamma\setminus \Gamma_2$.\end{lem}
\paragraph{}Let $\alpha,\beta \in \Gamma$. Write $\alpha \to \beta$ if there exists $l \geq 0$ such that $T^l(\alpha)=\beta$. We extend this relation to $\Z^+\Gamma$ by setting $\alpha \to \beta$ if there exists $\alpha_1,\ldots \alpha_r,\beta_1,\ldots \beta_r \in \Gamma$ such that $\alpha_i\to\beta_i$ for $i=1,\ldots r$ and $\alpha=\sum_i \alpha_i, \;\beta=\sum_i\beta_i$. It is easy to see that this relation is transitive, i.e if $\alpha \to \beta$ and $\beta \to \gamma$ then $\alpha \to \gamma$. Set 
$$\Z^+\Gamma_{\to\alpha}=\{\sigma \in \Z^+\Gamma, \sigma \to \alpha\},\qquad 
\Z^+\Gamma_{\alpha\to}=\{\sigma \in \Z^+\Gamma, \alpha\to \sigma\}.$$
Now let us consider the space
$$I_{T}=\bigoplus_{\beta \to \alpha}\big(U_q(\n_-)[-\alpha]q^{(\Z^+\Gamma_{\to\alpha})^\vee}\otimes U_q(\n_+)[\beta]q^{(-\Z^+\Gamma_{\beta\to})^\vee}\big)\subset U_q(\b_-) \otimes U_q(\b_+).$$
\begin{lem}The space $I_{T}$ is stable under the actions of $B \otimes 1$, $1 \otimes B^{-1}$ and $Ad(q^{\Omega_\h})$.\end{lem}
\noindent
\textit{Proof.} Note that the actions of $(B \otimes 1)$ and $(1 \otimes B^{-1})$ are well-defined on $I_{T}$ as $B(U_q(\n_-)[-\alpha])=0$ if $\alpha \not\in \Z^+\Gamma_1$ and $B^{-1}(U_q(\n_+)[\beta])=0$ if $\beta \not\in \Z^+\Gamma_2$. It is clear that $(B \otimes 1)I_{T}\subset I_{T}$ and $(1 \otimes B^{-1})I_{T}\subset I_{T}$. The last claim in the Lemma follows easily from the formula
$$Ad(q^{\Omega_\h})(u_\alpha \otimes v_{\beta})=u_\alpha q^{\beta^\vee}\otimes q^{-\alpha^\vee}v_\beta$$
if $u_\alpha \in U_q(\n_-)[-\alpha]$ and $v_\beta \in U_q(\n_+)[\beta]$.\qed
\paragraph{}Note that the Cayley transform $C_T: \h_0 \to \h_0$ is still well-defined in the Kac-Moody setting. Set $Z=\frac{1}{2}((1-C_T)\otimes 1)\Omega_{\h_0}$. Let $\overline{I}_{T}$ be the completion of $I_{T}$ with respect to the principal gradings in $U_q(\b_{\pm})$ and let $\overline{I}^*_{T}$ be the subspace consisting of elements of strictly negative degree in the first component and strictly positive degree in the second component. 
\begin{theo} There exists a unique element $\mathcal{J}^0_T(\lambda):\l^* \to (1+\overline{I}^*_{T})^\l$ such that
\begin{equation}\label{E:MABRR0}
\mathcal{R}^{21}q_1^{2\lambda}B_1\mathcal{J}^0_T(\lambda)=\mathcal{J}^0_T(\lambda)q_1^{2\lambda}q^{\Omega_\h}.
\end{equation}
Moreover $\mathcal{J}_T(\lambda):=\mathcal{J}^0_T(\lambda)q^Z$ satisfies the 2-cocycle relation
$$\mathcal{J}^{12,3}_T(\lambda)\mathcal{J}^{12}_T(\lambda+\frac{1}{2}h^{(3)})=\mathcal{J}^{1,23}_T(\lambda)\mathcal{J}^{23}_T(\lambda-\frac{1}{2}h^{(1)}).$$
\end{theo}
\noindent
\textit{Proof.} The first statement is proved exactly as in \cite{ESS}. We write $\mathcal{J}^0_T(\lambda)=1+\sum_{j \geq i}\mathcal{J}^{0,i}_{T}(\lambda)$ where $\mathcal{J}^{0,i}_{T}(\lambda)$ has degree $i$ in the first component. Then (\ref{E:MABRR0}) is equivalent to a system of equations labelled by $j \geq 1$
$$Ad(q^{\Omega_\h}q_1^{2\lambda})B_1 \mathcal{J}^{0,j}_T(\lambda)=\mathcal{J}^{0,j}_T(\lambda)+\ldots$$
where $\dots$ stands for terms involving $\mathcal{J}^{0,i}_T(\lambda)$ with $i<j$. But the operator $Ad(q^{\Omega_\h}q_1^{2\lambda})B_1 -1$ is invertible on $I_{T}^\l$ for generic $\lambda$ and $\mathcal{J}^{0,j}_{T}(\lambda)$ can be computed recursively.\\
\hbox to1em{\hfill}The second claim is proved as \cite{ESS}, Section 4. We consider the three components versions of (\ref{E:MABRR0})
\begin{equation}\label{E:KM1}
\mathcal{R}^{21}\mathcal{R}^{31}q_1^{2\lambda}B_1X^0_T(\lambda)=X^0_T(\lambda)q_{12}^{\Omega_\h}q_{13}^{\Omega_\h},
\end{equation}
\begin{equation}\label{E:KM2}
\mathcal{R}^{32}\mathcal{R}^{21}q_3^{-2\lambda}B_3^{-1}X^0_T(\lambda)=X^0_T(\lambda)q_{12}^{\Omega_\h}q_{13}^{\Omega_\h},
\end{equation}
acting on (a suitable completion of) the space
$$\bigoplus_{\alpha,\beta,\gamma}\big(U_q(\n_-)[-\alpha]q^{(\Z^+\Gamma_{\to \alpha})^\vee}\otimes U_q(\g)[\beta]\otimes U_q(\n_+)[\gamma]q^{(-\Z^+\Gamma_{\gamma\to})^\vee}\big)$$
where the sum runs over all triples $(\alpha,\beta,\gamma)$ such that $\beta$ can be written as $\beta=\beta^+-\beta^-$ where $\beta^++\gamma \to \beta^-+\alpha$. It is not difficult to show that $(\mathcal{J}^0_T(\lambda))^{1,23}\mathrm{Ad}\; q^Z_{1,23}(\mathcal{J}^0_T(\lambda+\frac{1}{2}h^{(1)}))^{23}$ and $(\mathcal{J}^0_T(\lambda))^{12,3}\mathrm{Ad}\; q^Z_{12,3}(\mathcal{J}^0_T(\lambda+\frac{1}{2}h^{(3)}))^{12}$ are two solutions of (\ref{E:KM1}) and (\ref{E:KM2}) with the same degree zero terms (in component $1$ or in component $3$). This implies that they are equal (see \cite{ESS}, Lemma 4.3).\qed
\paragraph{}Now let $V_1,\ldots,V_N$ be $U_q(\g)$-modules from the category $\mathcal{O}$. Define the renormalized twisted traces functions $F^T_{V_1,\ldots,V_N}(\lambda,\mu)$ in the same way as in Section 2. Note that all the operators $\mathcal{J}_T(\lambda)$, $\mathbb{R}_T(\lambda)$, $\mathbb{Q}_T(\lambda)$,... are well-defined on any module from the category $\mathcal{O}$ when considered as formal powers series in $q^{2(\lambda,\mu)}\C[[q^{-(\lambda,\alpha_i)}, q^{-(\mu,\alpha_i)}]]$, $\alpha_i \in \Gamma$. Operators $\mathcal{D}_W$ for affine algebras $\g$ are defined in some particular situation in \cite{E3}.

\begin{theo}\label{T:801}
The function $F^T_{V_1,\ldots, V_N}(\lambda,\mu)$ satisfies the following difference equation for all $j=1,\ldots, N$~:
\begin{equation}
F^T_{V_1,\ldots, V_N}(\lambda,\mu)=(D^T_j \otimes K^T_j) F_{V_1,\ldots, V_N}^T(\lambda,\mu)
\end{equation}
where $D^T_j$ and $ K^T_j$ are defined by (\ref{E:0pi21}) and (\ref{E:0pi22}).
\end{theo}
\begin{theo}\label{T:802}
Let $T$ be an automorphism of $\Gamma$. The functions $F^T_{V_1,\ldots, V_N}(\lambda,\mu)$ satisfy the following difference equation for each $j=1\ldots, N$ :
\begin{equation}
B_{V_j}B^*_{V_j^*}F_{V_1,\ldots, V_j,\ldots, V_N}^T(\lambda,\mu)=(D^{\vee,T}_j\otimes K^{\vee,T}_j) F^T_{V_1,\ldots, V_j^B,\ldots, V_N}(\lambda,\mu),
\end{equation}
where $D^{\vee,T}_j$ and $ K^{\vee,T}_j$ are defined by (\ref{E:71}).
\end{theo}
The above two theorems are proved in the same way as Theorems~\ref{T:02} and \ref{T:04} respectively.

\paragraph{}Similarly, let $W$ be an integrable highest weight $U_q(\g)$-module (resp. a $B$-invariant integrable highest weight $U_q(\g)$-module) and let $V_1,\ldots,V_N$ be $U_q(\g)$-modules from the category $\mathcal{O}$. 
\begin{theo}\label{T:803}
\begin{equation}
\mathcal{D}_W^T F^T_{V_1,\ldots, V_N}(\lambda,\mu)=\chi_W(q^{-2\mu})F^T_{V_1,\ldots, V_N}(\lambda,\mu),
\end{equation}
where $\mathcal{D}_W^T$ is defined by (\ref{E:0pi1}).
\end{theo}
\begin{theo}\label{T:804}
 Let $T$ be an automorphism of $\Gamma$. Then
\begin{equation}
\mathcal{D}^{\vee,T}_{W}F^T_{V_1,\ldots, V_N}(\lambda,\mu)=\mathrm{Tr}_{|W^{\h_0}}(q^{-2\lambda}B)F^T_{V_1,\ldots, V_N}(\lambda,\mu),
\end{equation}
where $\mathcal{D}^{\vee,T}_{W}$ is defined by (\ref{E:70}).
\end{theo}
The proof of the above two theorems is the same as in the finite-dimensional case.
\paragraph{Remark.} The integrability condition on the module $W$ is not essential.
\paragraph{}The classical limits of Theorems~\ref{T:801}-\ref{T:804} are analogous to the the corresponding classical limits of Theorems~\ref{T:01}-\ref{T:04} in Section 7.\\
\\ 
\centerline{\textbf{Acknowledgments}}
The first author was partially supported by the NSF grant DMS-9700477. The work of both authors was partly done when they were employed by the Clay Mathematics Institute as CMI Prize Fellows. O.S. would like to thank the MIT mathematics department for its hospitality.
\small{}
\noindent Pavel Etingof, MIT Mathematics Dept., 77 Massachusetts Ave.,
 CAMBRIDGE 02139 MA., USA\\
\texttt{etingof@math.mit.edu}\\
\noindent Olivier Schiffmann, MIT Mathematics Dept., 77 Massachusetts Ave.,
 CAMBRIDGE 02139 MA., USA\\
\texttt{schiffma@math.mit.edu}.
 \end{document}